\documentclass[12pt, letterpaper, twoside]{article}
\usepackage[utf8]{inputenc}
\usepackage[T1]{fontenc}
\usepackage[pdftex]{graphicx}
\usepackage[bookmarks=true]{hyperref}
\usepackage{xcolor}
\usepackage{indentfirst}
\usepackage{amsmath}
\usepackage{amsfonts}
\usepackage{amssymb}
\usepackage{array}
\usepackage{comment}
\usepackage{mathtools} 
\usepackage{bm}
\usepackage[width = 16cm, height = 22cm]{geometry}
\usepackage{siunitx}
\usepackage{setspace}
\usepackage{lineno}

\newcommand{\VRVE}{\ensuremath{V_0}}

\newcommand{\dd}[1]{\mathrm{\,d}\hspace{0.05em}#1}

\newcommand{\ndim}{\ensuremath{n_{\mathrm{dim}}}}

\usepackage[style=numeric,sorting=none,backend=biber, uniquename=false ,maxbibnames=30, maxcitenames=2]{biblatex}
\usepackage{empheq}

\newcommand{\footremember}[2]{%
\footnote{#2}
\newcounter{#1}
\setcounter{#1}{\value{footnote}}%
}
\newcommand{\footrecall}[1]{%
\footnotemark[\value{#1}]%
} 

\newcounter{JEremark}

\newcounter{HYremark}

\begin{document}

\pagestyle{plain}
\pagenumbering{arabic}

\title{Homogenization of discrete diffusion models by asymptotic expansion}
\author{Jan Eliáš\footremember{Brno}{Brno University of Technology, Faculty of Civil Engineering, Brno, Czechia}\footremember{email}{Corresponding author: jan.elias@vut.cz} \and Hao Yin\footremember{Northwestern}{Northwestern University, Department of Civil and Environmental Engineering, Evanston, IL USA} \and Gianluca Cusatis\footrecall{Northwestern}}
\date{}

\maketitle 

\section*{Abstract}
 
Diffusion behaviors of heterogeneous materials are of paramount importance in many engineering problems. Numerical models that take into account the internal structure of such materials are robust but computationally very expensive. This burden can be partially decreased by using discrete models, however even then the practical application is limited to relatively small material volumes. 

This paper formulates a homogenization scheme for discrete diffusion models. Asymptotic expansion homogenization is applied to distinguish between (i) the continuous macroscale description approximated by the standard finite element method and (ii) the fully resolved discrete mesoscale description in a local representative volume element (RVE) of material. Both transient and steady-state variants with nonlinear constitutive relations are discussed. In all the cases, the resulting discrete RVE problem becomes a simple linear steady-state problem that can be easily pre-computed. The scale separation provides a significant reduction of computational time allowing the solution of practical problems with a~negligible error introduced mainly by the finite element discretization at the macroscale. 

\section*{Key words}
homogenization, mass transport, diffusion, discrete model, concrete, Poisson's equation, quasi-brittle material

\section{Introduction}

Transport of mass, heat, or electric charge through heterogeneous solid material is an~important problem in a~number of practical applications. A large number of experimental and numerical works concerning heterogeneous materials are devoted to study the transport of moisture or water~\parencite{NguRah-19,BesMat-22}, heat~\parencite{JbeZha21,Sak21} or chlorides~\parencite{AlfAls19,WuWan-20}. Particularly in civil engineering, the transport behavior of concrete is of paramount importance as temperature, relative humidity, and/or chemical substances largely influence its mechanical behavior as well as durability~\parencite{BazJir18}.   

Concrete is a~heterogeneous material composed of mortar matrix and mineral aggregates. Its transport properties are, to a~large extent, dictated by the internal structure. Robust modeling approaches therefore involve explicit incorporation of the heterogeneity into the model. Unfortunately, such mesoscale models are extremely computationally demanding when real-size structures or structural members are being analyzed. The most efficient mesoscale model is arguably the discrete particle model, where the mesostructure is simplified into rigid units representing mineral grains and compliant contacts accounting for soft matrix. These kinds of models were originally developed for prediction of mechanical behavior of concrete and related materials, e.g.~\parencite{CusPel-11,CusMen-11,Eli16,AvaJir-21}. A~thorough description of the discrete modeling approach is provided by a~recent review paper by \textcite{BolEli-21}. The coupling of mechanics with mass transport phenomena is routinely implemented within the described discrete approach.

\textcite{BerBol04} developed an~elegant solution to the transport problem within the framework of discrete mechanical models where conduit elements coincide with the mechanical ones. This approach was later updated by the usage of a~dual lattice created as a~counterpart of the mechanical one~\parencite{NakSri-06,Gra09,GraBol16}. The dual approach provides an~easy coupling scheme capable to directly take into account crack openings in the transport problem because conduit elements run along the cracks. For example, the dual network is used to simulate hydraulic fracturing~\parencite{AsaPan-18,LiZho-18,UlvSun18,PhaSor-21}, deterioration due to the expansion of corrosion products~\parencite{FahWhe-17}, alkali-silica or other chemical reactions ~\parencite{AlnCus-13,AlnLuz-17,WanMen-19,TonWan20}, chloride diffusion~\parencite{WanUed11,SavLuk-14,TraPha-21}, or thermal spalling~\parencite{SheLi-20,SheZha-22}. The dual network can be enhanced by the presence of short fiber reinforcement~\parencite{PanKan-21}. But even for such efficient mesoscale models, the computational requirements become unbearable in many practical applications. Reduction of this burden is possible via computational homogenization. 

Homogenization is a~classical technique in analyzing heterogeneous materials~\parencite{BenLio-78,SanZao87}. It relies on two fundamental assumptions: (i) existence of representative volume element (RVE) or periodic unit cell, a~material volume containing complete information about the material internal structure and its properties; (ii) scale separation, an~assumption that the size of such RVE is much smaller than the size of the actual domain of the problem.
Computational homogenization~\parencite{SmiBre-98},  sometimes called FE$^2$, has been developed to overcome complications related to nonlinear constitutive models or complex material internal structures. Vast literature about computational homogenization in heterogeneous materials is available and it spans from the homogenization of mechanical properties~\parencite{GhoLee-95,FisQin-99,ChoSem-21}, heat conduction~\parencite{OzdBre-08,WasHeu-20}, transport phenomena~\parencite{BenSte20,SykKre-12}, flow~\parencite{KolPat-18} and diffusion~\parencite{PolLar-21}, to coupled processes~\parencite{SanFre-16}. 

Motivated by a~homogenized solution for discrete mechanical models \parencite{RezCus16,RezZho-17,LiRez-17,RezAln-19}, the asymptotic expansion homogenization is employed here to develop a~homogenization scheme for the discrete diffusion problem. Verification is performed with two different mesoscale flow lattice models. The first one has a~geometry of the dual transport network derived from Voronoi tessellation, the second model adopts the geometry from the Lattice-Discrete Particle Model (LDPM)~\parencite{CusPel-11,CusMen-11}. The main difference between those flow models is the geometrical relation between the normal $\mathbf{o}$ and contact vector $\mathbf{e}_{\lambda}$ for each conduit element (Fig.~\ref{fig:2Dsketch}b). The Voronoi tessellation renders them parallel, while these vectors generally differ when LDPM tessellation is used.

The derivation of the homogenization is performed specifically for mass transport with pressure being the primary field. It is however valid for a~large number of problems defined at continuum level by Poisson's equation at steady state or its transient extension. These include, for example, diffusion and heat or electric conduction.   

Several simple examples featuring linear or nonlinear steady-state and transient behavior are included at the end of the paper to verify the derived solution. There is also an application example showing analysis of coupled heat and moisture evolution during concrete dam curing and service life. 

\section{Discrete model of mass transport \label{sec:model}}

\begin{figure}[!tb]
	\centering
	\includegraphics[width=12cm]{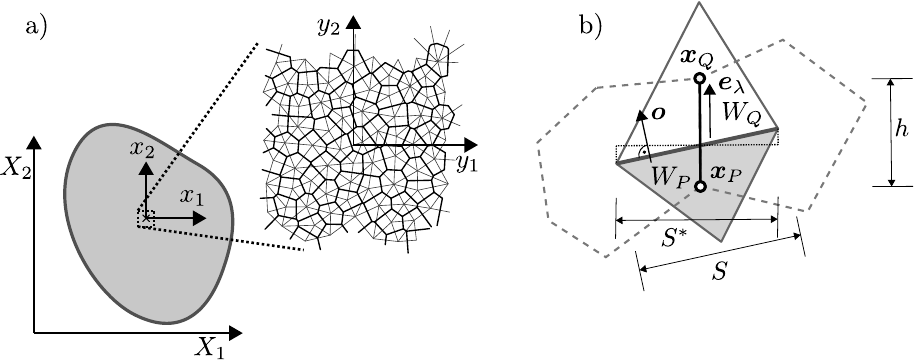}
	\caption{2D sketch of the homogenization setup: a) different reference systems, continuous macroscale problem and discrete periodic microscale problem; b) conduit element along the contact between two particles of the primary network. \label{fig:2Dsketch}}
\end{figure}

The primary variable of our reference problem in the discrete framework is pressure $p$. In contrast to the mechanics of the solid state, compressive pressure is considered positive. 
Each model node $P$ at coordinates $\mathbf{x}_P$ with pressure degree of freedom $p^P$ is surrounded by facets (forming a~simplex of the dual network called here a~control volume) that connect it to other nodes. The connection between nodes $P$ and $Q$ has some area $S$, length $h$ taken as Euclidean distance between the nodes, and direction (or contact vector) $\mathbf{e}_{\lambda}  = (\mathbf{x}_Q-\mathbf{x}_P)/h$. One also defines vector $\mathbf{x}_{PQ}=\mathbf{x}_Q-\mathbf{x}_P$. The volume of the simplex associated with node $P$ is denoted $W_P$. The situation is sketched in 2D in Fig.~\ref{fig:2Dsketch}b. 

A~discrete estimation of the pressure gradient in the normal direction between two nodes $P$ and $Q$ reads
\begin{align}
g = \nabla p \cdot \mathbf{e}_{\lambda} \approx \frac{p^Q-p^P}{h} \label{eq:presgrad}
\end{align} 
The constitutive equation relating the scalar flux density and the pressure gradient reads
\begin{align}
j = -\lambda(p_{\lambda})g \label{eq:LDPMconstFlow} \end{align}
with permeability coefficient $\lambda$ dependent on some weighted average pressure $p_{\lambda}$ of the two nodes $P$ and $Q$.

The planar facet of the conduit element has normal $\mathbf{o}$. If this vector is parallel to the contact vector ($\mathbf{o}||\mathbf{e}_{\lambda}$), the total flux through the element (from $P$ to $Q$, in the direction of $\mathbf{e}_{\lambda}$, outside from the control volume) can be estimated as $Sj$. This is precisely the case of models with geometry based on Voronoi or power/Laguerre tessellation~\parencite{Eli20}. Otherwise, as in the case of LDPM tessellation, one can project the true area $S$ into the plane defined by the contact vector $\mathbf{e}_{\lambda}$ and compute the total flux as
$S^{\star}j$ where $S^{\star} = |\mathbf{o}\cdot\mathbf{e}_{\lambda}| S$.

The mass balance equation is established for each control volume separately based on the conservation of mass. In the transient case, it reads  
\begin{align}
 W c(p) \dot{p} + \sum_{Q\in W} S^{\star} j =  W q(p) \label{eq:balance}
\end{align} 
for control volume $W$ formed by elements connecting node $P$ to surrounding nodes $Q$ (gray areas in Fig.~\ref{fig:2Dsketch}b). 
$c$ is the material capacity and $q$ is the source or sink term. The capacity, as well as the source term, is assumed to be dependent on pressure $p$. For the sake of simplicity, the balance equation omits indices $P$ and $Q$ identifying which flow element the quantities belong to. The source term $q$ is positive when the mass flows into the volume, flux density $j$ is positive when mass flows out.

The balance equations are assembled for all control volumes, appropriate boundary conditions are supplemented and the system is solved. Implicit integration schemes are used to integrate time in the transient analysis; in this study, the flow lattice based on LDPM tessellation uses the Backward Euler method while the model based on Voronoi tessellation uses the Generalized-$\alpha$ method~\parencite{CarKuh12}.

\section{Separation of scales}

Three spatial reference systems are considered: (i) global macroscopic reference system $\mathbf{X}$, (ii) local macroscopic (slow) reference system $\mathbf{x}$ unique for each RVE and (iii) local microscopic (fast) reference system $\mathbf{y}$ determining position within the RVE (see Fig.~\ref{fig:2Dsketch}a). The following scale separation relationship is assumed
\begin{align}
\mathbf{x} &= \eta \mathbf{y}  \label{eq:scalesep}
\end{align} 
with $\eta$ being the separation of scales constant with properties $0<\eta\ll1$. It defines the contrast between the macro and micro scales. 

All the geometric variables introduced in Sec.~\ref{sec:model} are defined in $\mathbf{x}$ reference system. They need to be transferred into the $\mathbf{y}$  reference system to reveal what power of $\eta$ they contain. These variables transformed into the $\mathbf{y}$ reference system are denoted by the $\tilde{\bullet}$ accent. The following transformation rules assign each variable an~appropriate power of $\eta$ reflecting power of distance unit involved 
\begin{align}
\mathbf{x}_{PQ} &= \eta{\mathbf{y}}_{PQ} & h &= \eta \tilde{h} & S^{\star} &= \eta^2 \tilde{S}^{\star} &
W &= \eta^3 \tilde{W} & \label{eq:scalingFlow}
\end{align}

The primary variable is now decomposed into the macroscopic (slow) component $p^{(0)}$ and microscopic (fast) component $p^{(1)}$ while terms with higher power of $\eta$ are omitted
\begin{align}
p = p^{(0)} + \eta p^{(1)} + \dots   \label{eq:decomposition}
\end{align} 
In accordance with the homogenization theory, it is assumed that all these terms are periodic over some RVE with periodic geometrical structure and the volumetric average of the fluctuation term $p^{(1)}$ over the RVE domain is zero.

From the viewpoint of the macroscopic spatial coordinate $\mathbf{X}$, neighboring nodes $P$ and $Q$ are close to each other. The value of pressure at point $Q$ can be approximated by Taylor expansion from point $P$  with respect to the global coordinate $\mathbf{X}$ \parencite{FisChe-07}.
\begin{align}
p^{QQ} = p(\mathbf{X}_Q, \mathbf{y}_Q) \approx p^{PQ} + \frac{\partial p^{PQ}}{\partial X_i} x_{i}^{PQ} + \frac{1}{2} \frac{\partial^2 p^{PQ}}{\partial X_i X_j} x^{PQ}_i x^{PQ}_j + \mathcal{O}(h^3) \label{eq:Taylorseries}
\end{align} 
A~pair of upper indices denotes the position in the $\mathbf{X}$ and $\mathbf{y}$ reference system, respectively. Terms with a~power of $h$ above 2 are omitted as $h$ is very small with respect to the global variable $\mathbf{X}$.

One substitutes Eq.~\eqref{eq:Taylorseries} into Eq.~\eqref{eq:presgrad} and obtains expression for pressure gradient in the macroscopic reference system $\mathbf{x}$
\begin{align}
g &= \frac{p^{QQ} - p^{PP}}{h} = \frac{1}{h}\left[p^{PQ} + \frac{\partial p^{PQ}}{\partial X_i} x_{i}^{PQ} + \frac{1}{2} \frac{\partial^2 p^{PQ}}{\partial X_i X_j} x^{PQ}_i x^{PQ}_j - p^{PP} \right] \label{eq:gbasic}
\end{align} 
Note that all the $\mathbf{X}$ coordinates are equal to $\mathbf{X}_P$, which allows us to remove one upper index and simplify the notation. Thus, only one upper index will be used hereinafter referring always to RVE at location $\mathbf{X}_P$: $p^{\bullet}=p^{P\bullet}=p(\mathbf{X}_P, \mathbf{y}_{\bullet})$. Furthermore, the upper index will indicate what $\eta$ power is involved. For example, $p^{0Q}=p^{(0)}(\mathbf{X}_P,\mathbf{y}_Q)$.

Next, equation~\eqref{eq:decomposition} is substituted into Eq.~\eqref{eq:gbasic}. All the length variables are transformed from the macroscopic system $\mathbf{x}$ to the microscopic reference system $\mathbf{y}$ according to Eq.~\eqref{eq:scalingFlow} and components of pressure gradient $g$ of the same $\eta$ power are collected. Three distinct terms in $g$ expression emerge (the fourth one with second spatial derivative of $p^{(1)}$ is omitted as negligible)
\begin{align}
g & = \eta^{-1} g^{(-1)} + g^{(0)} + \eta g^{(1)} + \dots
\end{align}
where
\begin{subequations} \label{eq:gterms}
\begin{align}
g^{(-1)} &=\frac{1}{\tilde{h}}\left[p^{0Q} - p^{0P}\right] \label{eq:g-1}\\
g^{(0)} &= \frac{1}{\tilde{h}}\left[p^{1Q} +  \frac{\partial p^{0Q}}{\partial X_i} y^{PQ}_i -  p^{1P} \right] \label{eq:g0}\\
g^{(1)} &= \frac{1}{\tilde{h}}\left[\frac{\partial p^{1Q}}{\partial X_i} y^{PQ}_i +\frac{1}{2} \frac{\partial^2 p^{0Q}}{\partial X_i X_j} y^{PQ}_i  y^{PQ}_j\right] \label{eq:g1} 
\end{align} 
\end{subequations}
For later use, it is useful to compute the gradient of $g^{(0)}$ with respect to the global variable $\mathbf{X}$
\begin{align}
\nabla_X g^{(0)} = \frac{\partial g^{(0)}}{\partial X_i} &= \frac{1}{\tilde{h}}\left[\frac{\partial \left(p^{1Q}-p^{1P}\right)}{\partial X_i} +  \frac{\partial^2 p^{0Q}}{\partial X_i X_j} y^{PQ}_j \right] \label{eq:gradj0}
\end{align} 
In the derivation of Eq.~\eqref{eq:gradj0}, $\partial y^{PQ}_i/\partial X_j=0$ because the spatial derivative of $\mathbf{y}_{PQ}$ with respect to the global coordinate $\mathbf{X}$ is zero due to the periodic RVE geometry.

The flux density is also decomposed into three contributions
\begin{align}
j & = -\lambda\left( p^{(0)}_{\lambda}+\eta p^{(1)}_{\lambda} \right) \left( \eta^{-1} g^{(-1)} +  g^{(0)} + \eta g^{(1)} + \dots\right) = \eta^{-1} j^{(-1)} + j^{(0)} + \eta j^{(1)} + \dots \label{eq:jdecomp}
\end{align}
The permeability coefficient $\lambda$ is expanded using Taylor series at $p=p^{(0)}_{\lambda}$
\begin{align}
j & = -\left(\lambda\left( p^{(0)}_{\lambda}\right) + \eta  \frac{\partial \lambda(p)}{\partial p}\bigg|_{p=p^{(0)}_{\lambda}} p^{(1)}_{\lambda} + \dots \right)\left( \eta^{-1} g^{(-1)} +  g^{(0)} + \eta g^{(1)} + \dots \right) \label{eq:jdecompexpanded}
\end{align}
The first flux term $j^{(-1)}$ must be, because of the scale separation, equal to 
\begin{align}
j^{(-1)} = -\lambda\left(p^{(0)}_{\lambda} \right) g^{(-1)}
\end{align}
It will be shown later that $j^{(-1)}$ is always zero and therefore $g^{(-1)}$ is always zero as well
\begin{align}
j^{(-1)} = g^{(-1)} = 0 \label{eq:j-1} 
\end{align}
The second term, after dropping the zero 
term $g^{(-1)}$, becomes the actual constitutive function for flow defined in Eq.~\eqref{eq:LDPMconstFlow} with permeability coefficient  dependent only on the macroscopic (slow) pressure $p^{(0)}$
\begin{align}
j^{(0)} & = -\lambda\left(p^{(0)}_{\lambda}\right) g^{(0)} \label{eq:j0}
\end{align}
The third term, $j^{(1)}$, reads from Eq.~\eqref{eq:jdecompexpanded}
\begin{align}
j^{(1)} & = -\lambda\left(p^{(0)}_{\lambda}\right)g^{(1)} - \frac{\partial \lambda(p)}{\partial p}\bigg|_{p=p^{(0)}_{\lambda}} p^{(1)}_{\lambda}g^{(0)} \label{eq:jexpanded}
\end{align}
which can be rewritten using Eq.~\eqref{eq:j0}
\begin{align}
j^{(1)} & = \frac{\partial j^{(0)}}{\partial g^{(0)}} g^{(1)} +\frac{\partial j^{(0)}}{\partial p^{(0)}_{\lambda}} p^{(1)}_{\lambda} \label{eq:j1}
\end{align}

Finally, the Taylor series expansions of the capacity and source terms around $p^{(0)}$ read
\begin{subequations} \label{eq:cqexpanded}
\begin{align}
c\left(p^{(0)}+\eta p^{(1)}+\dots\right) & = c\left(p^{(0)}\right) + \eta \frac{\partial c}{\partial p}\bigg|_{p=p^{(0)}} p^{(1)} + \dots \label{eq:cexpanded} \\
q\left(p^{(0)}+\eta p^{(1)}+\dots\right) & = q\left(p^{(0)}\right) + \eta \frac{\partial q}{\partial p}\bigg|_{p=p^{(0)}} p^{(1)} + \dots \label{eq:qexpanded}
\end{align}
\end{subequations}
For convenience, symbols $c_0 = c\left(p^{(0)}\right)$ and $q_0 = q\left(p^{(0)}\right)$ represent parts of the capacity and source term dependent solely on the macroscopic primary variable $p^{(0)}$.

\section{Balance equations}
Assuming that density $\rho_w$, capacity $c$, permeability coefficient $\lambda$ and source terms $q$ are of order $\mathcal{O}(\eta^0)$, one can write down the balance equation~\eqref{eq:balance} with all variables transformed to $\mathbf{y}$ reference system. After dividing by $\eta^3$, the balance equation reads
\begin{align}
\tilde{W} c \dot{p} + \frac{1}{\eta} \sum_{Q\in W} \tilde{S}^{\star} j  = \tilde{W} q 
\end{align} 
The balance equation, after substituting Eqs.~\eqref{eq:decomposition}, \eqref{eq:jdecomp} and \eqref{eq:cqexpanded}, can be decomposed into three separate equations collecting the terms with appropriate power of $\eta$. The terms with $\eta$  power above 0 are omitted.

\subsection{Constant pressure $p^{(0)}$ within RVE ($\eta^{-2}$)} 
The terms with negative second power of $\eta$, transferred back to $\mathbf{x}$ reference system  and multiplied by $\eta^3$, read
\begin{align}
\eta^{-1}\sum_{Q\in W} S^{\star} j^{(-1)} = 0 \label{eq:balance-2}
\end{align} 
Assembling such equation for each control volume, one obtains system of equations with trivial solution 
\begin{align}
j^{(-1)} = 0 \label{eq:solution-2}
\end{align} 
Recalling that $j^{(-1)}$ depends solely on $g^{(-1)}$ via some admissible constitutive relation, it follows that  (using Eq.~\ref{eq:g-1})
\begin{align}
g^{(-1)} = \frac{1}{\tilde{h}}\left[p^{0Q} - p^{0P}\right] = 0
\end{align} 
and therefore 
\begin{align}
p^{0P} = p^{0Q} \label{eq:p0const}
\end{align} 
This shows that the macroscopic primary field $p^{(0)}$ does not depend on the local variable $\mathbf{y}$ and proves the previously assumed Eq.~\eqref{eq:j-1}. The $\mathbf{y}$ coordinate upper index is therefore dropped hereinafter in expressions belonging to $p^{(0)}: p^{0Q}=p^{0P}=p^{(0)}$. Also the average pressure within the element can be simplified: $p^{(0)}_{\lambda} =p^{(0)}$. 

The system of equation is linear in the case on a~linear constitutive model. One can then easily prove that the solution~\eqref{eq:solution-2}-\eqref{eq:p0const} is unique. Even though some other solutions to Eq.~\eqref{eq:balance-2} may theoretically exist in the nonlinear regime, the reported trivial solution will always be assumed. Such solution would indeed occur because it would always be energetically more favorable configuration.

\subsection{RVE level ($\eta^{-1}$) \label{sec:RVEproblem}} 
The negative first power of $\eta$ collects the following terms (already transformed to the $\mathbf{x}$ reference system and multiplied by $\eta^3$)
\begin{align}
\sum_{Q\in W}S^{\star}j^{(0)} = 0 \label{eq:RVEbalance}
\end{align}
with $j^{(0)}$ depending on $g^{(0)}$ and $p^{(0)}$ via Eq.~\eqref{eq:j0}. The pressure gradient $g^{(0)}$ defined in Eq.~\eqref{eq:g0} is transformed to the $\mathbf{x}$ reference system
\begin{align}
g^{(0)} &= \frac{\eta}{h}\left[p^{1Q} - p^{1P} \right] - \hat{g}  \label{eq:gradRVE}
\end{align} 
with $\hat{g}$ being the eigen term representing the projected gradient of the macroscopic solution
\begin{align}
\hat{g}= -\nabla_X p^{(0)} \cdot \mathbf{e}_{\lambda} = -\mathbf{a} \cdot \mathbf{e}_{\lambda} \label{eq:eigenstrainRVE}
\end{align}
where $\mathbf{a}=\nabla_X p^{(0)}$. Since $g^{(0)}$ depends on $\eta p^{(1)}$ via Eq.~\eqref{eq:gradRVE}, Eq.~\eqref{eq:RVEbalance} provides the balance equations for the periodic RVE level where the degrees of freedom are the pressures $\eta p^{(1)}$. Note that there is no transient term in Eq.~\eqref{eq:RVEbalance}, the RVE level is a~steady-state boundary value problem. Also, according to the Eq.~\eqref{eq:j0}, the permeability coefficient $\lambda$ is dependent solely on the macroscale (slow) primary field $p^{(0)}$. Therefore, the problem that was originally transient and nonlinear becomes steady-state and linear at the RVE level. The ``load'' is provided by the eigen term in Eq.~\eqref{eq:eigenstrainRVE} and comes from the gradient of the macroscopic pressure $p^{(0)}$.

The $\eta p^{(1)}$ primary field must be periodic, therefore periodic boundary conditions are enforced. Furthermore, one additional Dirichlet boundary condition, determining some absolute value of pressure, is needed (otherwise, only differences in the pressure are involved). The zero volumetric average of the fluctuation field $p^{(1)}$ shall be prescribed, $\langle p^{(1)}\rangle = 0$, where the the weighted volumetric average reads 
\begin{align}
\langle\bullet\rangle = \frac{1}{\VRVE} \sum_{P\in \VRVE} W \bullet \label{eq:volumAver}
\end{align} 
However, such a~condition is not straightforward to implement. Moreover, if it is implemented using linear constraints, the originally sparse matrix of the RVE system becomes full and memory requirements as well as computational time to solve the RVE problem rapidly grow with RVE size. 

It is therefore advised to randomly pick a~node and set its pressure to some prescribed value (for example, zero). One then obtains exactly the same fluxes and pressure gradients, only the primary pressure field $\eta p^{(1)}$ is shifted by an~unknown constant. Optionally if one is interested in the true values of the microscale (fast) pressure, the zero volumetric average can be enforced during post-processing.

\subsection{Macroscopic level ($\eta^0$)} 
Collecting the terms with $\eta^0$ results in the following equation (already transformed to the $\mathbf{x}$ reference system and multiplied by $\eta^3$)
\begin{align}
\eta \sum_{Q\in W} S^{\star} j^{(1)} =  W q_0 - W c_0 \dot{p}^{(0)}
\label{eq:macrobalance}
\end{align} 
Summing all the terms over the whole RVE volume  at coordinates $\mathbf{X}_P$, dividing the result by RVE volume $\VRVE$ and substituting Eq.~\eqref{eq:j1} into Eq.~\eqref{eq:macrobalance} provides
\begin{align}
\frac{\eta}{\VRVE}\sum_{P\in \VRVE}\sum_{Q\in W} S^{\star} \left( \frac{\partial j^{(0)}} {\partial g^{(0)}}g^{(1)} + \frac{\partial j^{(0)}}{\partial p^{(0)}} p^{(1)}_{\lambda} \right) = \frac{\sum_{P\in \VRVE} W q_0}{\VRVE} - \frac{\sum_{P\in \VRVE} W c_0 }{\VRVE} \dot{p}^{(0)}  = 0  \label{eq:macrobalancesummed}
\end{align} 
Each term on the left-hand side appears in the summation twice, once connecting node $P$ with node $Q$ and once connecting $Q$ with $P$. The contact vectors have opposite directions in these cases ($\mathbf{e}_{PQ}=-\mathbf{e}_{QP}$), but their product is the same as two negative signs cancel each other. Also $j^{(0)}$ and $g^{(0)}$ have opposite signs while $\partial j^{(0)}/\partial g^{(0)}$, $p^{(0)}$ and $p^{(1)}_{\lambda}$ are identical for these two terms. 

The pressure gradients $g^{(1)}$ belonging to elements $PQ$ and $QP$ are now summed (using Eqs.~\ref{eq:g1} and \ref{eq:gradj0}) while the average pressures $p^{(1)}_{\lambda}$ are subtracted 
\begin{subequations} 
\begin{align}
^{PQ}g^{(1)} + ^{QP\!\!}g^{(1)} &= \frac{1}{\tilde{h}} \left[\frac{\partial \left(p^{1Q} - p^{1P}\right)}{\partial X_i}  y^{PQ}_i + \frac{\partial^2 p^{(0)}}{\partial X_i X_j} y^{PQ}_i y^{PQ}_j \right] = \frac{\partial g^{(0)}}{\partial X_i} y^{PQ}_i\\
^{PQ}p^{(1)}_{\lambda} - ^{QP\!\!}p^{(1)}_{\lambda} &= 0
\end{align} 
\end{subequations}
Substituting this result back to Eq.~\eqref{eq:macrobalancesummed} along with the definition of the weighted volumetric average~\eqref{eq:volumAver} 
and using the chain rule,
one obtains the following  balance equation at the macroscopic level ($\mathbf{y}_{PQ}$ is replaced by $\mathbf{x}_{PQ}/\eta$)
\begin{align}
\frac{1}{\VRVE}\sum_{e\in \VRVE} S^{\star}\frac{\partial j^{(0)}}{\partial X_i} x^{PQ}_i  =  \langle q_0 \rangle - \langle c_0 \rangle \rho_w \dot{p}^{(0)}  \label{eq:macrobalanceaveraged}
\end{align} 
Note that by summing the contributions of two identical terms with the opposite direction, the double summation changed into a~single summation over all conduit elements $e$ in the RVE.

The chain rule for divergence provides
\begin{align}
\nabla_X\cdot\left( j^{(0)}S^{\star}\mathbf{x}_{PQ}\right) &= S^{\star}\nabla_X j^{(0)} \cdot \mathbf{x}_{PQ}
\end{align} 
because the divergence of the geometric quantities with respect to the global variable $\mathbf{X}$ is zero due to the periodic RVE geometry: $\nabla_X\cdot (S^{\star} \mathbf{x}_{PQ})=0$. This allows us to define the macroscopic flux vector  
\begin{align}
\mathbf{f} = \frac{1}{\VRVE}\sum_{e\in\VRVE}  h S^{\star} j^{(0)} \mathbf{e}_{\lambda} \label{eq:macroflux}
\end{align} 
and to simplify Eq.~\eqref{eq:macrobalanceaveraged} to be the standard continuous diffusion equation
\begin{align}
\langle c_0 \rangle \rho_w \dot{p}^{(0)} + \nabla_X\cdot \mathbf{f} =  \langle q_0 \rangle \label{eq:macrofinal}
\end{align} 
where the degrees of freedom of the macroscopic problem are pressures $p^{(0)}$. If the same material is used within the whole RVE domain, both capacity and source term can be simplified, because they depend only on the constant primary variable $p^{(0)}$  (Eq.~\ref{eq:cqexpanded}): $\langle c_0 \rangle=c_0$ and $\langle q_0 \rangle=q_0$. The macroscopic problem described by Eq.~\eqref{eq:macrofinal} can be easily solved by the Finite Element Method. Transformation into the weak form is straightforward. The material routine evaluated at each integration point is replaced by the RVE problem described in Sec.~\ref{sec:RVEproblem}.

\section{Simplifications under additional assumptions \label{sec:discussion}}

The derived equations describe a~general solution suitable for any functional form of permeability coefficient $\lambda(p)$ as long as it depends only on the pressure magnitude. If there is a~dependence on other variables, they need to be added to the Taylor expansion of $\lambda$ in Eq.~\eqref{eq:jexpanded}. The same is valid for capacity and source term in Eqs.~\eqref{eq:cqexpanded}. Also, any spatial randomness in the material parameters at the level of the RVE can be directly incorporated.

It is worth investigating some special cases listed below.

\subsection{Pre-computed conductivity tensor and fast solution \label{sec:precomputed}}
One can build an~extremely efficient solution approach based on a~special form of the constitutive model, which is however widely applied. If the nonlinear dependence of the permeability coefficient on pressure $p^{(0)}$ takes a~form of multiplication of some material constant by a~relative nonlinear effect, $\lambda = \lambda_0 \kappa_{\mathrm{r}}\left(p^{(0)}\right)$, the RVE solution can be easily pre-computed. The parameter $\lambda_0$, the intrinsic permeability coefficient, is a~user defined material constant that can vary within the RVE, while the relative permeability $\kappa_{\mathrm{r}}\left(p^{(0)}\right)$ provides the nonlinear dependence identical to all flow elements. 

Let us combine the RVE equations~\eqref{eq:j0}, \eqref{eq:RVEbalance}, \eqref{eq:gradRVE} and \eqref{eq:eigenstrainRVE} into one 
\begin{align}
\sum_{Q\in W}S^{\star} \lambda\left(p^{(0)}\right) \left( \frac{\eta}{h}\left[p^{1Q} - p^{1P} \right] + \mathbf{a} \cdot \mathbf{e}_{\lambda} \right) = 0 \label{eq:RVEcombined}
\end{align}
where $\mathbf{a}=\nabla_X p^{(0)}$ is the macroscopic pressure gradient. By substituting the multiplicative decomposition of permeability coefficient $\lambda = \lambda_0 \kappa_{\mathrm{r}}\left(p^{(0)}\right)$, and since the pressure $p^{(0)}$ is constant within the RVE, the relative permeability  $\kappa_{\mathrm{r}}\left(p^{(0)}\right)$ must be constant as well and can be easily eliminated from the problem
\begin{align}
\sum_{Q\in W}S^{\star} \lambda_0 \left( \frac{\eta}{h}\left[p^{1Q} - p^{1P} \right] + \mathbf{a} \cdot \mathbf{e}_{\lambda} \right) = 0
\end{align}
The linear problem can now be solved in advance. Applying artificial pressure gradient repetitively in a~form of standard Cartesian basis unit vectors $\mathbf{i}$, one can pre-compute a set of solutions $p^{(1)}_i$ and fluxes $j^{(0)}_i$. The real solution for the actual pressure gradient $\mathbf{a}$ with entries $a_i$ is then simply a~linear combination of the pre-computed solution vectors; the flux density must be multiplied by the nonlinear term
\begin{align}
p^{(1)} &= \sum_{i=1}^{\ndim} a_i p^{(1)}_i    & j^{(0)} &= \kappa_{\mathrm{r}}\left(p^{(0)}\right) \sum_{i=1}^{\ndim} a_i j^{(0)}_i   \label{eq:fastfields}
\end{align}
where $n_{\mathrm{dim}}$ is the number of dimensions of the problem. Substituting this decomposition into Eq.~\eqref{eq:macroflux} yields
\begin{align}
\mathbf{f} &= \kappa_{\mathrm{r}}\left(p^{(0)}\right) \sum_{i=1}^{\ndim} a_i \mathbf{f}_i & \mathrm{with}\quad \mathbf{f}_i &= \frac{1}{\VRVE}\sum_{e\in\VRVE}  h S^{\star} j^{(0)}_i \mathbf{e}_{\lambda} \label{eq:precompf}
\end{align}
The pre-computed macroscopic flux vectors $\mathbf{f}_i$ can be row- or column-wise collected into a~second order symmetric conductivity tensor $\bm{\Lambda}$ of size $\ndim$, and one can add the negative sign for convenience. Because this tensor is identical for all RVEs and all time steps, it is enough to evaluate it only once at the beginning of the simulation and only for one RVE. It can be then used to easily compute the true nonlinear macroscopic flux in the RVEs
\begin{align}
\mathbf{f} &= -\kappa_{\mathrm{r}}\left(p^{(0)}\right) \bm{\Lambda} \cdot \mathbf{a} \label{eq:collectf}
\end{align}
This approach effectively transfers the nonlinearity to the macroscale by assuming that the nonlinear effects can be expressed as a~multiplication of the basic material permeability coefficient with spatially identical nonlinear function. 

\subsection{Constant permeability coefficient with parallel contact and normal vector \label{sec:constperm}}

The sum in the combined RVE equation~\eqref{eq:RVEcombined} is split and the eigen term is moved to the right-hand side
\begin{align}
\sum_{Q\in W}\eta\lambda\left(p^{(0)}\right)\frac{S^{\star}}{h}\left[p^{1Q} -  p^{1P} \right]  = -\mathbf{a}\cdot\sum_{Q\in W} \lambda\left(p^{(0)}\right) S^{\star} \mathbf{e}_{\lambda} \label{eq:RVERHS}
\end{align}
This shows that the projection of the macroscopic pressure can be introduced by assigning an~external source at each RVE node equal to the right-hand side of Eq.~\eqref{eq:RVERHS} instead of the eigen pressure gradient $\hat{g}$. However, the flux density $j^{(0)}$ needed to evaluate the macroscopic flux $\mathbf{f}$ in Eq.~\eqref{eq:macroflux} still depends on the eigen term via the pressure gradient $g^{(0)}$ defined in Eq.~\eqref{eq:gradRVE}.

If the same material is used everywhere within the RVE,  the $\lambda$ parameter is constant and can be eliminated from Eq.~\eqref{eq:RVERHS}. Furthermore, if the normal and contact vectors are parallel ($\mathbf{o}||\mathbf{e}_\lambda$), the right-hand side of Eq.~\eqref{eq:RVERHS} vanishes. Parallelism  $\mathbf{o}||\mathbf{e}_{\lambda}$, which is ensured when the model geometry is established by Voronoi or power/Laguerre tessellation, implies that there is no distortion of contact areas by the projection: $S^{\star}=S$. The right-hand side of Eq.~\eqref{eq:RVERHS} becomes 
\begin{align}
-\mathbf{a}\cdot\sum_{Q\in W} S \mathbf{e}_{\lambda}
\end{align}
The summation of area $S$ times outward normal $\mathbf{e}_{\lambda}$ is performed over enclosed surface $\Gamma$ of the control volume. Using the divergence theorem, one can show that it is always zero
\begin{align}
\sum_{Q\in W} S 
\mathbf{e}_{\lambda} = \int_{\Gamma} e^{\lambda}_i \dd{\Gamma} = \int_{\Gamma} \mathbf{i} \cdot \mathbf{e}_{\lambda} \dd{\Gamma} = \int_{W} \nabla\cdot \mathbf{i} \dd{W}  =\bm{0}  \label{eq:divergencetheorem1}
\end{align}
Therefore, for Voronoi or power tessellation and spatially constant permeability coefficient, the right-hand side of Eq.~\eqref{eq:RVERHS} disappears and the solution yields $\eta p^{(1)}=0$ everywhere. This exactly corresponds to the case of a homogeneous, isotropic material without any internal material textures or heterogeneities. 

The macroscopic flux vector then becomes 
\begin{align}
\mathbf{f} = \frac{1}{\VRVE} \lambda\left(p^{(0)}\right)  \sum_{e\in\VRVE}  h S \hat{g} \mathbf{e}_{\lambda} = - \frac{1}{\VRVE} \lambda\left(p^{(0)}\right) \mathbf{a} \cdot \sum_{e\in\VRVE}  h S \mathbf{e}_{\lambda} \otimes \mathbf{e}_{\lambda} \label{eq:fluxessimple1}
\end{align} 
To further modify this result, one can define the vector $\mathbf{r}$ as a~vector connecting the governing node (where the pressure degree of freedom is defined) with the centroid of the element facet. For each transport element, there are two such vectors, $\mathbf{r}_P$ and $\mathbf{r}_Q$, connecting the same centroid with two different governing nodes.  Projection of these vectors into the normal direction yields positive or negative distances $h_P=\mathbf{r}_P\cdot\mathbf{e}_{\lambda}$ and $h_Q=-\mathbf{r}_Q\cdot\mathbf{e}_{\lambda}$ from the nodes to the element contact plane. The vector $h\mathbf{e}_{\lambda}$ can be then rewritten as 
\begin{align}
h\mathbf{e}_{\lambda} = \mathbf{r}_P - \mathbf{r}_Q
\end{align}
and subsequently the whole summation in Eq.~\eqref{eq:fluxessimple1} reads 
\begin{align}
\sum_{e\in\VRVE}  h S \mathbf{e}_{\lambda} \otimes \mathbf{e}_{\lambda} = \sum_{e\in\VRVE} S\mathbf{r}_P\otimes \mathbf{e}_{\lambda} - S\mathbf{r}_Q\otimes \mathbf{e}_{\lambda}
\end{align} 
The summation can be reorganized by visiting each control volume $P$ only once and summing immediately all the contributions from all neighbors $Q$. Since the normal now always point in the outward direction, all the signs are positive.
\begin{align}
\sum_{e\in\VRVE}  h S \mathbf{e}_{\lambda} \otimes \mathbf{e}_{\lambda} = \sum_{P\in\VRVE}\sum_{Q\in{W}} S\mathbf{r}\otimes \mathbf{e}_{\lambda} \label{eq:helpsum}
\end{align} 
Finally, such summation over the enclosed surface can be, with the help of the divergence theorem, calculated analytically
\begin{align}
\sum_{Q\in W} S \mathbf{r}\otimes\mathbf{e}_{\lambda} = \int_{\Gamma} \mathbf{r} \otimes \mathbf{e}_{\lambda} \dd{\Gamma} = \int_{\Gamma} r_i \mathbf{j} \cdot \mathbf{e}_{\lambda} \dd{\Gamma} = \int_{W} \nabla \cdot (r_i \mathbf{j})\dd{W}  = \int_{W} \frac{\partial r_i}{\partial X_j}\dd{W} = \bm{1}W\label{eq:divergencetheorem2}
\end{align} 
where $\bm{1}$ denotes the second order identity tensor (or Kronecker delta). The second summation over all control volumes in Eq.~\eqref{eq:helpsum} yields the total sum to be $\bm{1} \VRVE$. Thus the Eq.~\eqref{eq:fluxessimple1} becomes
\begin{align}
\mathbf{f} = -\frac{1}{\VRVE}\lambda\left(p^{(0)}\right) \mathbf{a} \cdot \bm{1} \VRVE = -\lambda\left(p^{(0)}\right) \mathbf{a}
\end{align} 
For constant permeability coefficient and parallel normal and contact vectors, the slow solution $p^{(0)}$ is the only meaningful component, the whole RVE response (that requires the solution of a~nonlinear steady-state problem in a~general case) collapses into a~multiplication of the macroscopic pressure gradient $\mathbf{a}$ by microscopic user-defined permeability coefficient $\lambda\left(p^{(0)}\right)$. Again, this is the expected result for homogeneous, isotropic material internal structure.

\section{Verification of flow model with Voronoi tessellation}

This section verifies the derived equations by comparing the \emph{homogenized} solution with the \emph{full} solution utilizing the original discrete model in the whole domain. The calculations adopt a~discrete model with geometry provided by the Voronoi tessellation on randomly placed nuclei with minimal mutual distance $l_{\min}=10$ mm. When the RVE structure is generated, the periodic distance~\parencite{EliVor16} is used to ensure RVE periodicity. The following nonlinear model of permeability coefficient, originally developed by \textcite{Gen80}, is adopted from \parencite{GraBol16}
\begin{align}
\lambda = \frac{\rho_w \kappa_0}{\mu} \kappa_{\mathrm{r}}(z)
\end{align}
where $\mu$ is the viscosity of the liquid and $\kappa_0$ is the intrinsic permeability. The relative permeability $\kappa_{\mathrm{r}}$ is  dependent on saturation $z$ as
\begin{align}
\kappa_{\mathrm{r}}(z)  =\sqrt{z}\left[ 1 - \left(1-\sqrt[m]{z}\right)^2\right]^2 
\end{align}
and the saturation is estimated based on the pressure 
\begin{align}
z  =\left( 1 + \sqrt[1-m]{p/\alpha}\right)^{-m}
\end{align}
using two additional parameters $m$ and $\alpha$. Since this constitutive model obeys the assumption from Sec.~\ref{sec:precomputed}, the pre-computed solution was used in all the analyses discussed hereinafter. That means that all the analyses with the \emph{homogenized} model start with a~calculation of RVE fast fields $p_i$ and conductivity tensor $\bm{\Lambda}$ on a~single RVE assuming linearity of the constitutive model. The evaluation of flux vectors and conductivity tensors at integration points during the real solution of the problem then reduces to simple and extremely fast evaluation of Eqs.~\eqref{eq:fastfields}-\eqref{eq:collectf} depending only on macroscopic pressure gradient $\mathbf{a}$ and macroscopic pressure $p^{(0)}$.

For the numerical examples, the model parameters for the water diffusion through concrete are chosen to be the following: water dynamic viscosity $\mu=8.9\times10^{-4}\,$Pa$\cdot$s, water density $10^3\,$kg/m$^3$, intrinsic permeability $\kappa_0= 5\times10^{-18}$\,m$^{2}$~\parencite{ClaEsh-03}, $m=0.5$, $\alpha=1$\,MPa~\parencite{GraBol16}. Finally, the capacity parameter $c$ is set to $1.64\times10^{-5}\,$s$^2/$m$^2$~\parencite{LiZho-18}. 

As demonstrated in Sec.~\ref{sec:constperm}, the fast solution is for the case of spatially constant permeability coefficient and geometry based on Voronoi tessellation trivial and can be computed analytically. In order to avoid such simplicity, the intrinsic permeability parameter was randomized. Its values were drawn independently for each conduit element from a lognormal distribution with the mean value equal to $\kappa_0$ and a coefficient of variation equal to 0.2.

\subsection{Steady-state linear analysis}
The first example consists of a~simple two-dimensional prism of depth 0.3\,m and length 1.2\,m (virtual thickness is 1\,m) with pressure prescribed at both of its ends. The left end pressure is set to zero while the right end pressure linearly increases in time from 0 to 1\,MPa. The boundary conditions on the top and bottom surfaces prescribe zero vertical flux. 

\begin{figure}[!bt]
	\centering
	\includegraphics[width=12.5cm]{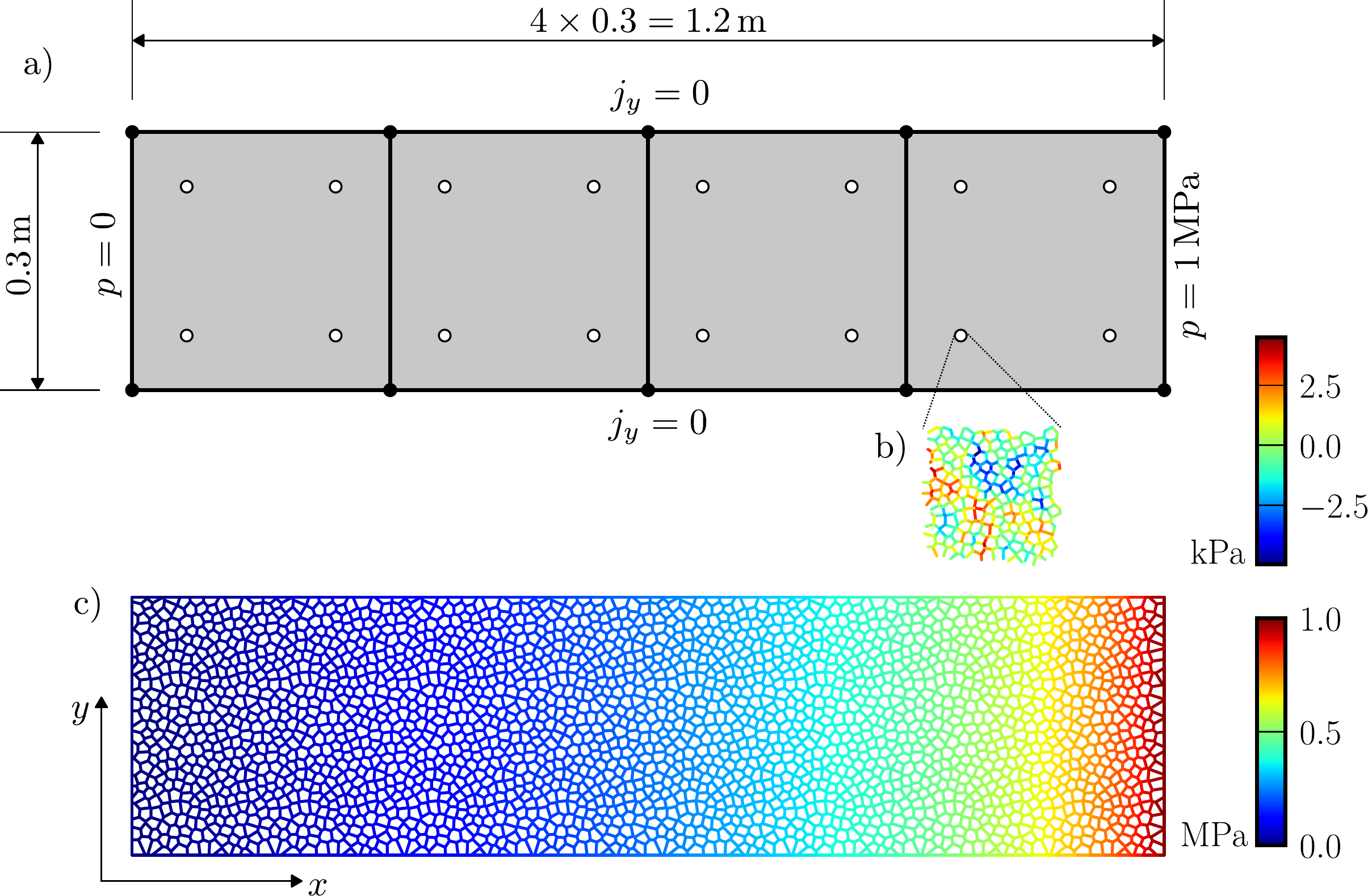}
	\caption{a) Dimensions, boundary conditions and division into finite elements, b) fast primary field $\eta p^{(1)}$ in the last time step in the periodic RVE attached to one of the integration points, c) \emph{full} solution in the last time step of steady-state nonlinear simulation.  \label{fig:bar}}
\end{figure}

The problem is first computed using the linear constitutive model where $\kappa_{\mathrm{r}}(z)=1$. The slow solution is approximated using 4-node bilinear isoparametric finite elements of size $0.3\times0.3$\,m (see Fig.~\ref{fig:bar}a). A~periodic RVE unit of size $0.15\times0.15$\,m is attached to each integration point (Fig.~\ref{fig:bar}b). The presented setup therefore uses the RVE size equal to the volume associated with the integration point. Such correspondence of volumes is extremely ineffective, the homogenized model shall be used with RVE size much smaller than the integration point volume. Nevertheless, it was chosen here to compare the \emph{homogenized} solution with the \emph{full} solution where the discrete model geometry is tiled from the actual RVEs (Fig.~\ref{fig:bar}c) including the spatial randomness of the permeability coefficient. Unfortunately the \emph{full} discrete structure slightly differs because geometry on the boundary of the \emph{full} prism must be adjusted~\parencite{Eli17} to exactly represent the prism surfaces. 

\begin{figure}[!tb]
	\centering
	\includegraphics[width=10cm]{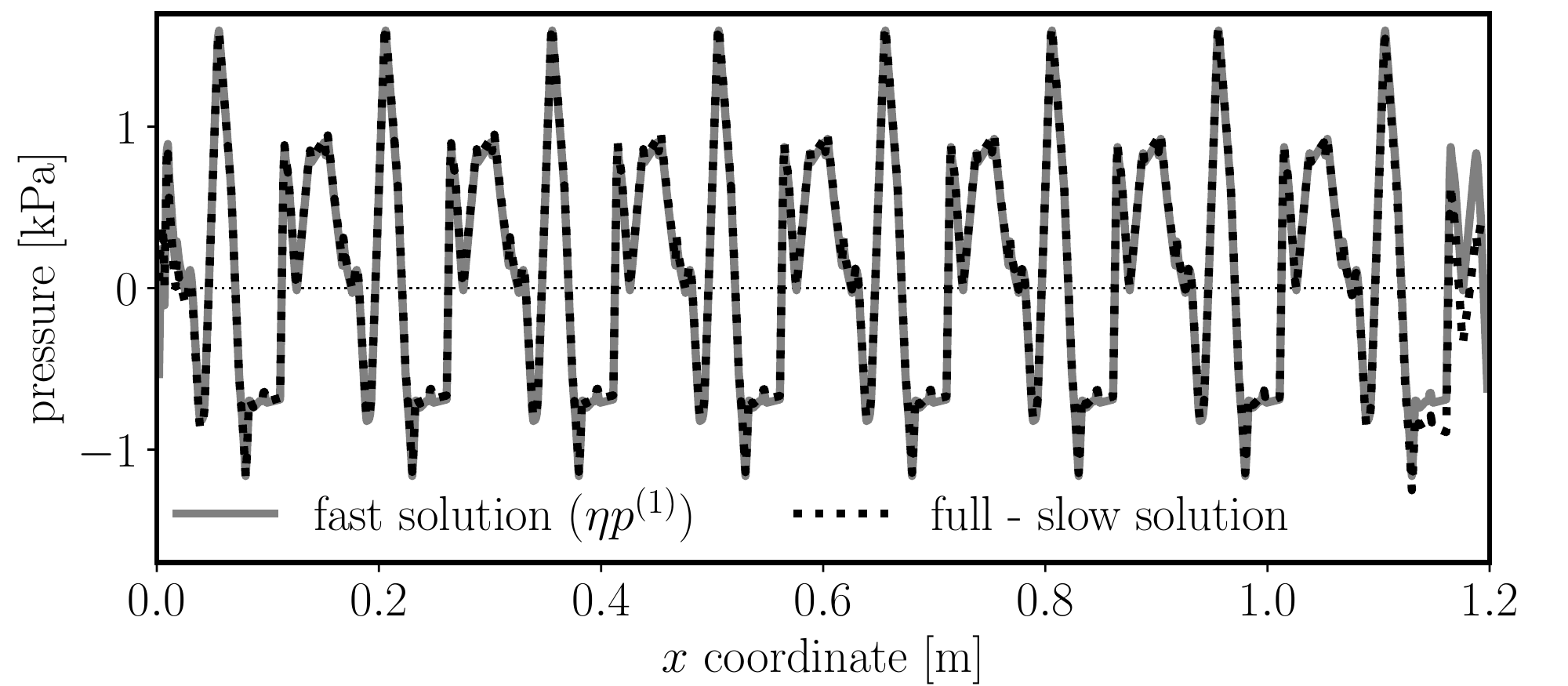}
	\caption{Fast pressure field $\eta p^{(1)}$ compared to the difference between the \emph{full} model pressure and slow pressure $p^{(0)}$ along the central line of the prism in the linear steady-state analysis.  \label{fig:profiles_linear}}
\end{figure}

Results of the steady-state analysis of the \emph{full} and \emph{homogenized} model are reported at the end time when pressure on the right-hand side reaches 1\,MPa. The total flux measured either at the left- or right-hand side of the prism equals 118.10\,g/day for the \emph{homogenized} and 118.02\,g/day for the \emph{full} solution, respectively. The pressure values from the RVE solution (i.e. $\eta p^{(1)}$ values) are interpolated to the central line of the prism (at $y=0.15$\,m), the slow solution $p^{(0)}$ along the same line is directly computed from the shape functions of the finite elements, and the \emph{full} solution is interpolated onto the same line as well. Figure~\ref{fig:profiles_linear} compares the $\eta p^{(1)}$ solution with the difference between the \emph{full} and $p^{(0)}$ solutions. The fast field oscillates in an~exact periodic manner as it is composed of repeating solutions at RVEs at macroscale integration points. The difference between the full and slow solutions exhibits almost identical oscillating pattern since the \emph{full} model geometrical structure is composed of repetitively inserted RVE structures. The correspondence is excellent except in the boundary regions where the internal structure of the RVE and \emph{full} prism differs because the tessellation is affected by the presence of boundaries~\parencite{Eli17}.

\begin{figure}[!tb]
	\centering
	\includegraphics[width=13cm]{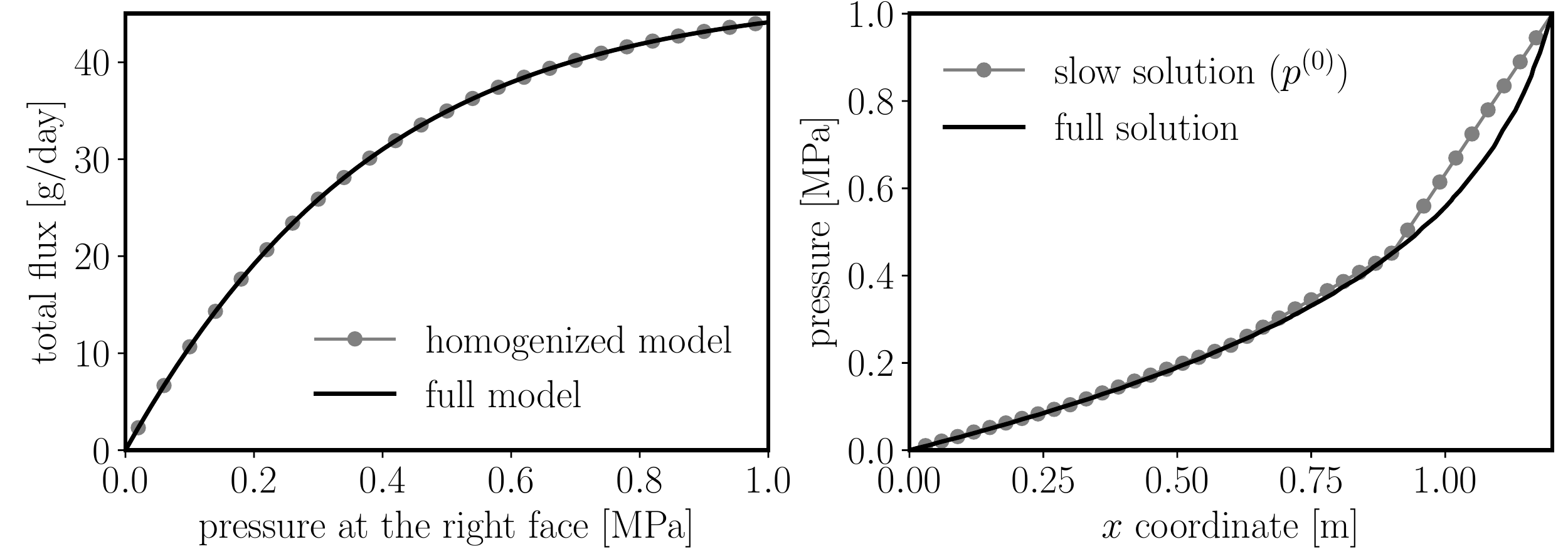}
	\caption{Verification of the homogenization method for nonlinear steady-state analysis: a) total flux at different times computed using the \emph{full} and \emph{homogenized} models; b) slow pressure field $p^{(0)}$ compared to the \emph{full} model pressure field along the central line of the prism. \label{fig:results_nonlinear}}
\end{figure}

\subsection{Steady-state nonlinear analysis}
The same setup is used to verify the nonlinear solution. The increase of the right-hand side boundary pressure to 1\,MPa is divided into 50 time steps. History of the total flux through the prism in time is shown in Fig.~\ref{fig:results_nonlinear}a. One can see that the correspondence between the \emph{full} and \emph{homogenized} solution is again excellent. There is however a~difference in the pressure profile along the central line. Figure~\ref{fig:results_nonlinear}b compares the slow solution $p^{(0)}$ to the \emph{full} solution and reveals that four linear finite element does not provide enough freedom to represent the highly nonlinear pressure profile along the prism with sufficient accuracy. 
The pressure increases towards the right-hand side and so the nonlinearity of the pressure profile. The element located at the far-right therefore exhibits the largest differences, because the linear shape functions cannot represent the correct curved profile. The errors are caused solely by the poor approximation of the macroscopic pressure field $p^{(0)}$. Using finer mesh or elements with higher polynomial approximations, the $p^{(0)}$ pressure would eventually match the curve in Fig.~\ref{fig:results_nonlinear}b, but the one-to-one match of the internal structure would be lost. Such convergence is shown in the next example verifying the transient results.

\begin{figure}[!tb]
	\centering
	\includegraphics[width=13cm]{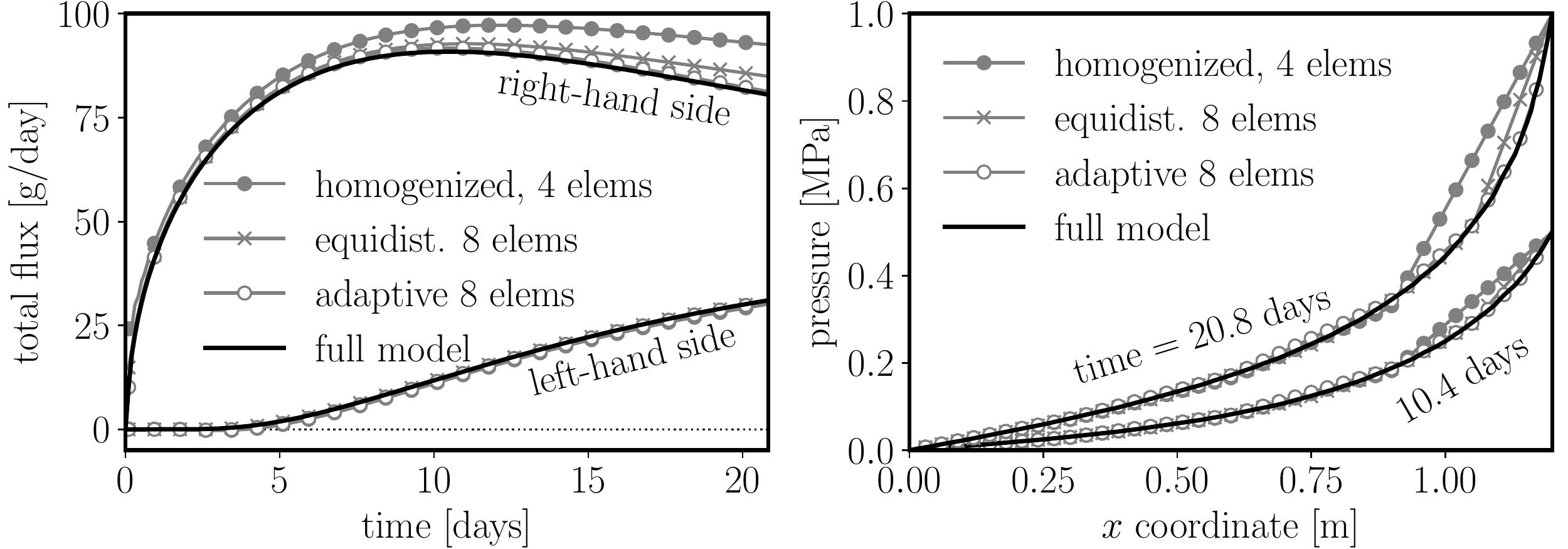}
	\caption{Total flux (left) and pressure field (right) along the central line computed in nonlinear transient analysis. Only macroscopic (slow) pressures $p^{(0)}$ are reported for the homogenized models.   \label{fig:results_transient}}
\end{figure}

\begin{figure}[!tb]
	\centering
	\includegraphics[width=13cm]{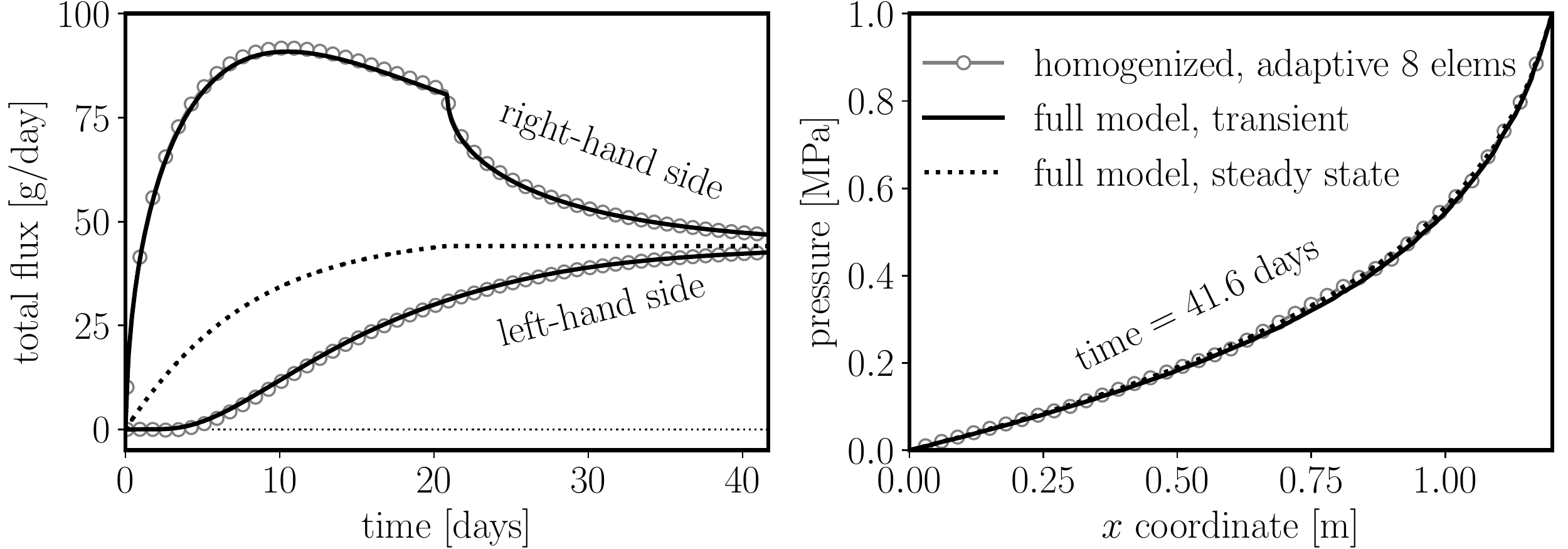}
	\caption{Total flux (left) and pressure field (right) computed in nonlinear transient analysis extended by a~period of transition towards steady state.  \label{fig:results_transient_to_steady}}
\end{figure}

\subsection{Transient analysis}
This section employs the same two-dimensional prism for the third time to verify the nonlinear transient analysis. The right-hand side pressure increases linearly in time from 0 to 1\,MPa in 200 time steps of length 9000\,s. The \emph{full} model is compared to the \emph{homogenized} model with different macroscopic primary field approximations. The first homogenized model builds the slow solution based on 4 isoparametric finite elements as sketched in Fig.~\ref{fig:bar}a. The second model improves it by using 8 elements of the same width. The third model, called here ``adaptive'', employs 8 elements as well but they have different widths --  proceeding from the left-hand side, there are 3 elements of width 0.3\,m, and then 5 elements of variable widths $w$, $w/2$, $w/4$, $w/8$ and $w/16$, with $w=24/155$\,m.  All of the homogenized models use the same RVE and therefore also the same pre-computed effective conductivity tensor based on linear steady-state analysis.

Fluxes measured at the left- and right-hand sides of the prism are plotted in Fig.~\ref{fig:results_transient} as well as pressure profile along the prism's center-line. The four element \emph{homogenized} model shows the largest error as its macroscopic pressure approximation does not have enough freedom to represent the true solution. Increasing the number of elements brings substantial improvement. Even better results are delivered by the adaptive model, where macroscopic degrees of freedom are shifted to the region where pressure profile nonlinearity is the strongest -- to the right-hand side. Then the macroscopic pressure field is finally capable to reasonably match the \emph{full} pressure profile. Consequently, the fluxes from the \emph{homogenized} model with the adaptive element width correspond to the fluxes from the \emph{full} model.

We let the \emph{full} and the adaptive \emph{homogenized} simulations run for another 200 time steps with constant pressure on the right-hand side of the prism equal to 1\,MPa. Both transient models gradually evolve towards the steady-state as shown in Fig.~\ref{fig:results_transient_to_steady}. Fluxes in time and pressure profiles at the final step are compared to nonlinear steady-state solutions from the previous section. Both transient solutions clearly converge to the correct steady-state.

\subsection{Effect of RVE size \label{sec:RVEsize}}
\begin{figure}[!b]
	\centering
	\includegraphics[width=10cm]{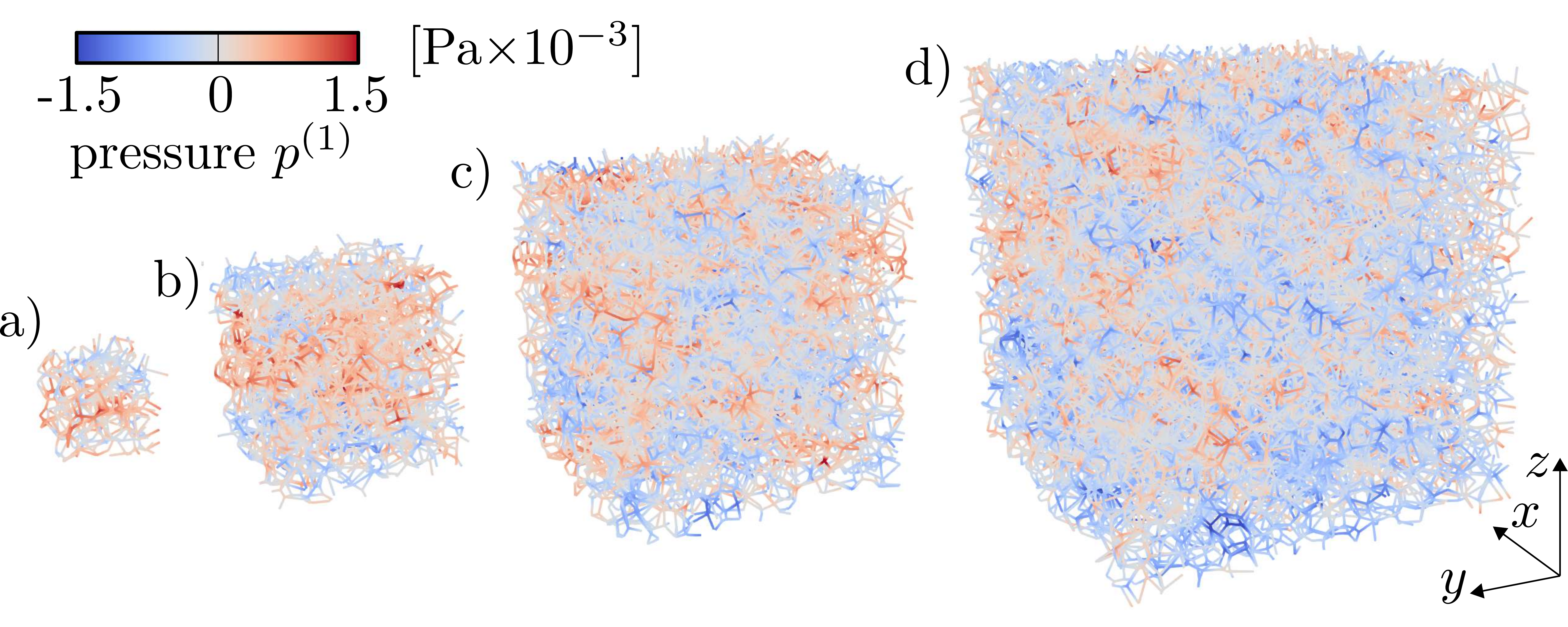}
	\caption{Pressure field $p^{(1)}$ in discrete cubical RVEs loaded by unit macroscopic pressure gradient in $z$ direction, $\nabla p^{(0)}=(0,\,0,\,1)$\,Pa/m. The RVE size and approximate number of degrees of freedom (DoF) are  a) 50\,mm, 500 DoFs, b) 100\,mm, 4000 DoFs, c) 150\,mm, 13500 DoFs and d) 200\,mm, 32000 DoFs.  \label{fig:cubes}}
\end{figure}

The size of the RVE should be as small as possible to save computational resources (especially when the RVE response is not pre-computed) but large enough to contain sufficient information about the average material response. To investigate the optimal RVE size, three-dimensional cubical RVEs of size 50, 100, 150 and 200\,mm were generated with the same Voronoi method ($l_{\min}=10$\,mm), the same material parameters and the same spatial randomness as before (Fig.~\ref{fig:cubes}). The whole RVE response can be, according to Sec.~\ref{sec:precomputed}, expressed by the 3$\times$3 symmetric tensor $\bm{\Lambda}$. Since the geometrical structure is isotropic, we expect that all diagonal entries should converge with size to identical values and all off-diagonal entries should vanish.

\begin{figure}[!tb]
	\centering
	\includegraphics[width=15cm]{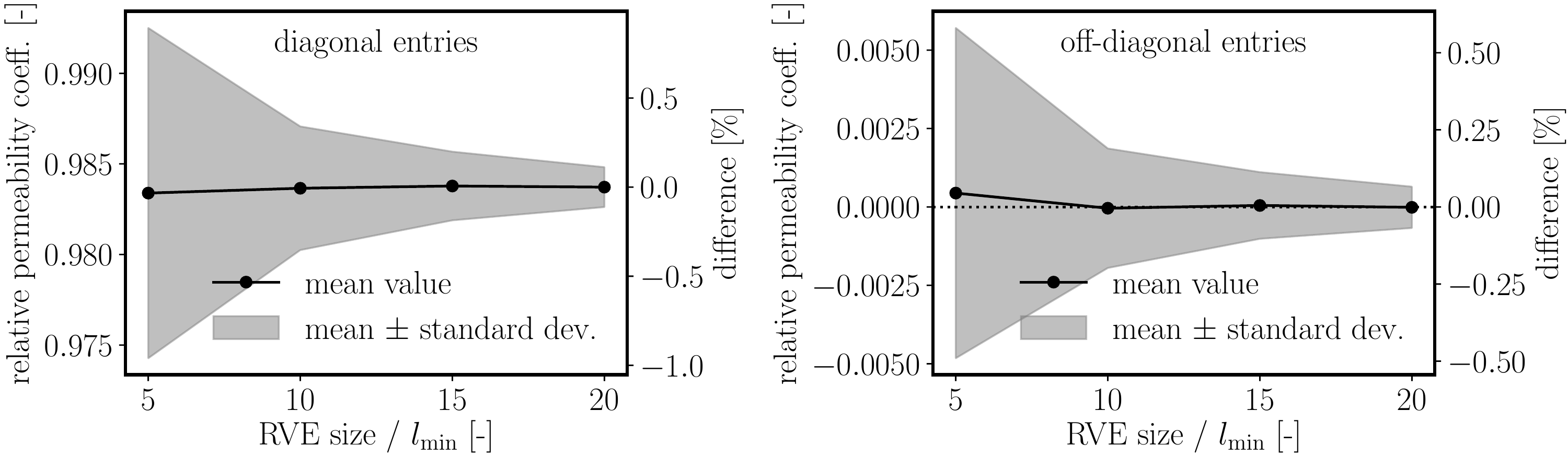}
	\caption{Diagonal (left) and off-diagonal (right) entries of the effective permeability tensor computed on Voronoi RVEs of different sizes.  \label{fig:RVE_size}}
\end{figure}

For each size 10 different internal structures were generated,  and for each internal structure 10 different random variants were considered. Each RVE variant resulted in 3 diagonal and 3 off-diagonal entries. Therefore in total, we computed 300 estimates of the diagonal and off-diagonal entries. The averages and standard deviations are plotted in Fig.~\ref{fig:RVE_size}. The left-hand side vertical axis shows flux (in or perpendicularly to the direction of the macroscopic pressure gradient, respectively) normalized by the mean permeability coefficient at the discrete element level and the magnitude of the pressure gradient. The right-hand side vertical axis shows relative difference using the average diagonal entry for RVE size 200\,mm as 100\%.

The averaged values are very similar irrespective of the RVE size. With a~larger number of realizations, these differences would diminish further. Also the average off-diagonal entries are close to zero as theoretically required. The optimum RVE size can be derived from the standard deviation, i.e., from the variability of the RVE behavior. One can choose how much variability is acceptable and select the corresponding RVE size. However, the relative differences are always below 1\% and therefore even the smallest RVE size seems to be sufficient for most practical applications.

Interestingly, the mean value of the diagonal entry does not converge to the intrinsic permeability coefficient $\lambda_0$, which is in our case 
$\lambda_0 = \rho_w \kappa_0  / \mu = 5.618\times10^{-12}\,$s (or 1 in relative units), but to some lower value around $5.526\times10^{-12}\,$s (or 0.984, respectively). This is due to the applied spatial randomness. If all the elements are connected in parallel, the effective permeability coefficient would be exactly equal to the mean of the applied distribution (which is $\lambda_0$). On the other extreme when all the elements are connected in series, the reciprocal value of the effective permeability coefficient would be equal to the mean of the reciprocal random permeability coefficients at the elements (i.e. harmonic mean). In our case of the lognormal distribution with a~coefficient of variation 20\%, the harmonic mean is numerically estimated to be $5.402\times10^{-12}\,$s (or 0.962 in relative units). Since the real RVE can be viewed as a~combination of parallel and serial connectivity, the effective permeability coefficient is found between those two extremes.

\section{Verification of flow model with LDPM tessellation \label{sec:verificationLDPM}}

The second flow lattice subjected to homogenization adopts the geometrical features of the Lattice Discrete Particle Model (LDPM)~\parencite{CusPel-11}. LDPM reflects concrete mesostructure by generating non-overlapping spherical particles of various diameters representing coarse aggregate particles. The maximum particle size is denoted as $d_{\max}$. The LDPM tessellation does not preserve parallelism between $\mathbf{o}$ and $\mathbf{e}_{\lambda}$ and therefore the projected elemental area differs from the true one: $S^{\star}\neq S$. As a~consequence, the simplification from Sec.~\ref{sec:constperm} does not hold and the spatial randomness introduced for the Voronoi model is now omitted.

The coupled Hygro-Thermo-Chemical (HTC)~\parencite{LuzCus09a,LuzCus09b} constitutive model will be used in this section. It simulates curing of cementitious material considering the chemical processes due cement hydration and silica fume reaction producing heat and consuming moisture along with simultaneous transient moisture and heat transport. The HTC model is briefly described in~\ref{sec:appendix} and features humidity and temperature dependent source, conductivity and capacity terms.  
The asymptotic expansion homogenization limits the coupling between temperature and humidity on the macroscopic (slow) components and therefore all the coupling effects occurs at the macroscale only. Also, the simplification introduced in Sec.~\ref{sec:precomputed} holds because the nonlinear parts of the constitutive model are spatially constant. This allows using the pre-computed RVE response and greatly reduces the computational time. 

\subsection{Effect of RVE size}
\begin{figure}[!b]
	\centering
	\includegraphics[width=15cm]{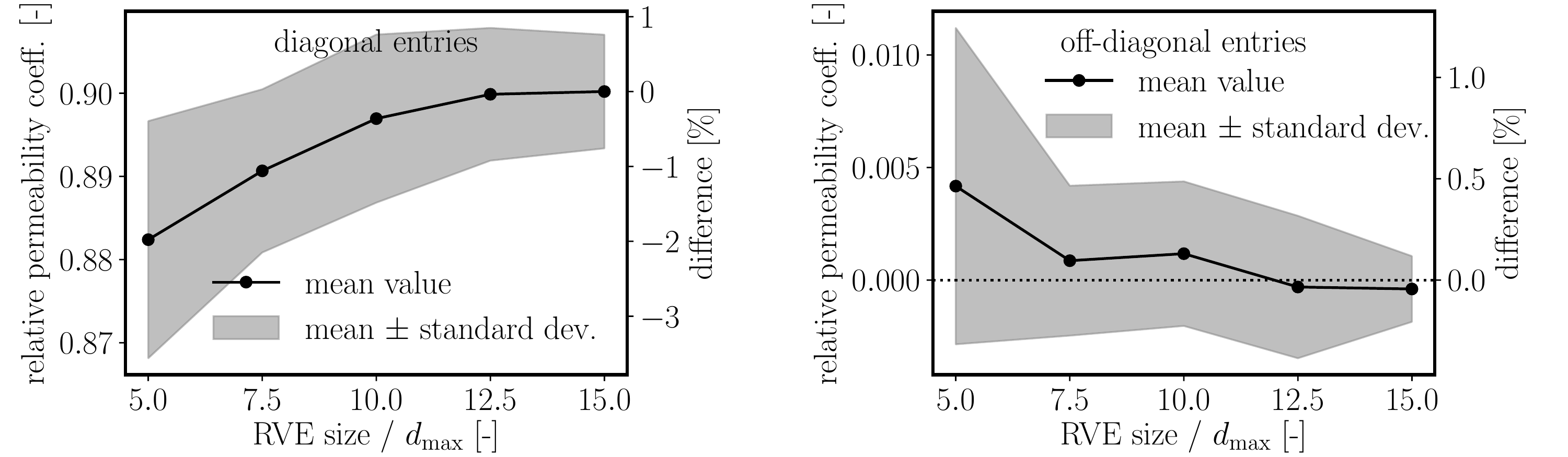}
	\caption{Diagonal (left) and off-diagonal (right) entries of the effective permeability tensor computed on RVEs of different sizes for flow lattice based on LDPM tessellation.  \label{fig:RVE_size_LDPM}}
\end{figure}

Before comparing the \emph{full} and \emph{homogenized} models, similar study as the one in Sec.~\ref{sec:RVEsize} of the effect of RVE size is performed here. Periodic RVE of sizes in range 100--300\,mm are generated using $d_{\max}=20$\,mm. The constitutive model used is a~simple linear one with an~arbitrary permeability coefficient, since the interest here is the linear and normalized response only. 10 different RVEs for each size were generated and analyzed, the average values and standard deviations of entries at and outside the diagonal of the permeability tensor are reported in Fig.~\ref{fig:RVE_size_LDPM}. The behavior is similar to what was observed for the Voronoi model. The off-diagonal entries vanishes with size as well as the variability of the diagonal entries.

There is, however, an~initial increase in the average value of the diagonal entries with RVE size which stabilizes for a~sufficient material volume. This increase is attributed to the boundary layer (wall effect) of a~similar type as described in Ref.~\parencite{Eli17} for mechanical models. The Voronoi model does not feature such an effect because it disappears when $\mathbf{o}$ and $\mathbf{e}_{\lambda}$ are parallel. Moreover, the Voronoi RVE was built without an~exact representation of the planar RVE boundaries. Instead, the boundaries followed surfaces of Delaunay simplices generated in a~random, periodic manner.

Note that the flow lattice with LDPM geometry does not reach a~relative effective permeability coefficient equal to 1 as well, even though there is no spatial variability considered. The reduced value of about 0.90 is due to projection of the area $S\rightarrow S^{\star}$ into the direction of the contact vector $\mathbf{e}_{\lambda}$. In effect, the effective permeability coefficient decreases below the permeability coefficient of the constitutive models at individual conduit elements.

\subsection{Comparison of full and homogenized models}

The comparison between \emph{full} and \emph{homogenized} solutions is now performed for a three-dimensional flow model with LDPM tessellation. The example employs the HTC model from \ref{sec:appendix} and it should approximately mimic the curing of a~concrete dam. The size of the concrete block is artificially small, allowing also the \emph{full} model to deliver results in a~reasonable time. The depth is 0.5\,m, the width at the bottom is 0.2\,m and the thickness is 0.2\,m, the maximum aggregate diameter is $d_{\max}=20$\,mm. The bottom surface, as well as the front and back surface are treated as sealed, no out-of-plane humidity or temperature flux is allowed. The temperature on the left- as well as on the right-hand and top surface is constant, $T_{\mathrm{w}}=T_{\mathrm{a}}=20^{\circ}$C. The relative humidity on the left-hand surface is set to $H_{\mathrm{w}}=1$ while the air relative humidity $H_{\mathrm{a}}$ on the right-hand and top surfaces is kept at 1 for 5 days and then linearly decreases within one day to 0.5 to represent demolding. From day 6, $H_{\mathrm{a}}=0.5$. The initial temperature and relative humidity in the small dam are $20^{\circ}$C and 1, respectively. The total simulation time is 1 year (52 weeks).  

\begin{figure}[!tb]
	\centering
	\includegraphics[width=\textwidth]{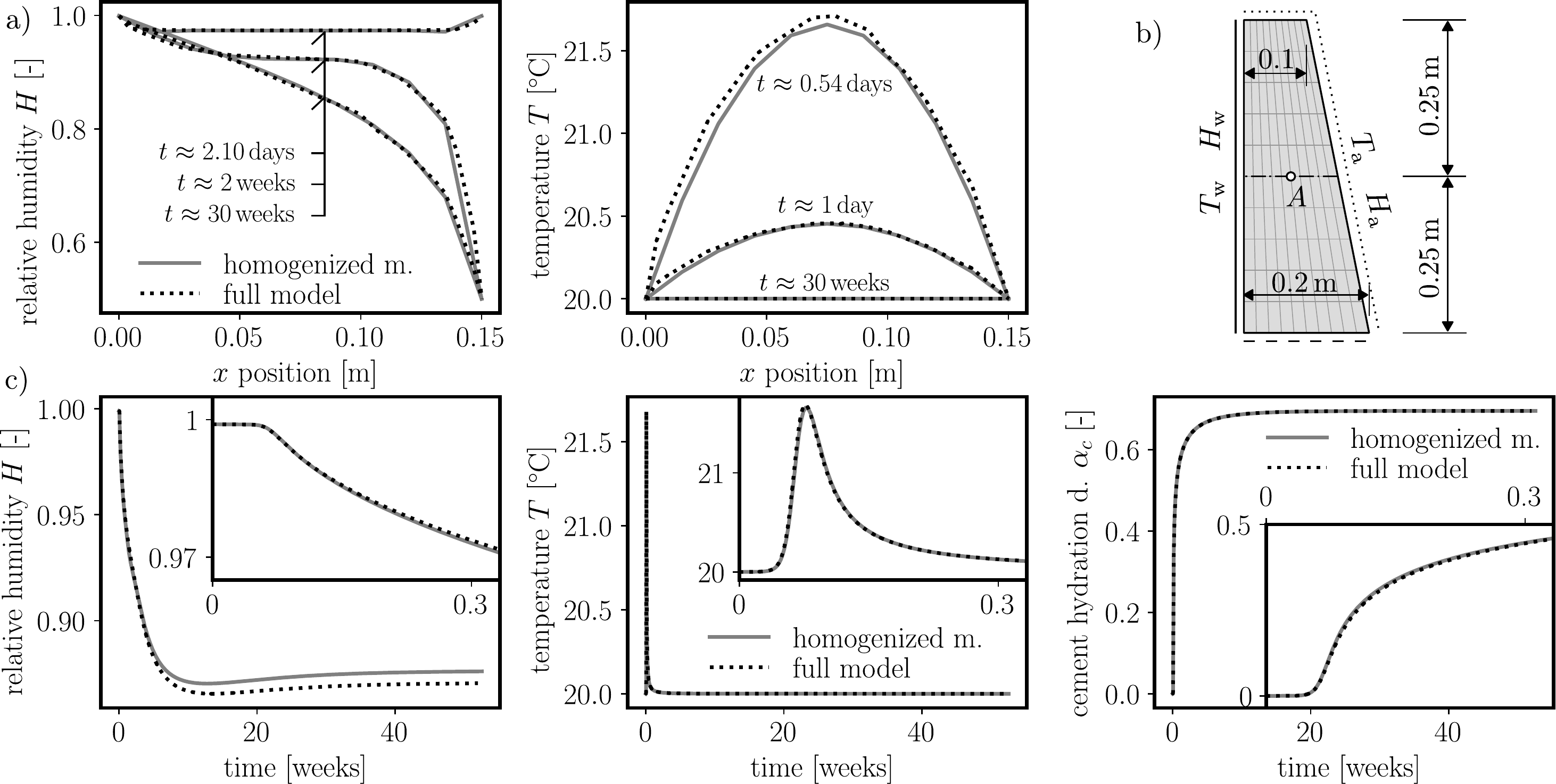}
	\caption{a) Humidity and temperature along the mid-depth horizontal line at different times; b) dimensions and boundary conditions; c) evolution of humidity, temperature and cement hydration degree at the central point $A$.  \label{fig:small_dam}}
\end{figure}

\begin{figure}[!tb]
	\centering
	\includegraphics[width=4.5cm]{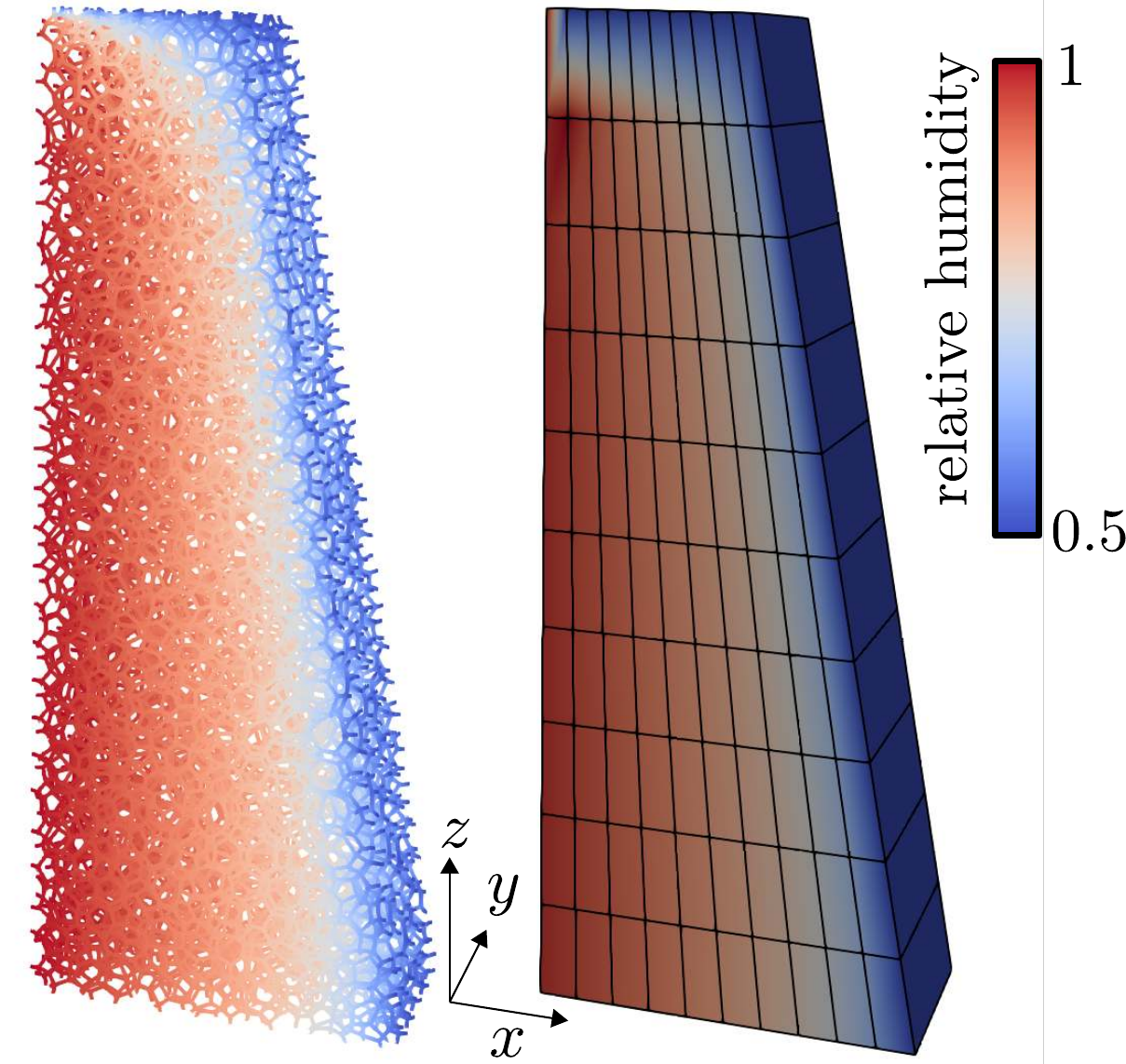}
	\caption{Spatial distribution of humidity in the small dam after one year: \emph{full} and \emph{homogenized} model. \label{fig:small_dam_paraview}}
\end{figure}

Figure~\ref{fig:small_dam}a shows relative humidity and temperature profiles along the mid-depth horizontal line (Fig.~\ref{fig:small_dam}b, dash-and-dot line) at different time instants. The \emph{homogenized} model shows just the macroscopic solution resulting from shape functions of finite elements while the \emph{full} model curves are calculated as averages over the nodal values in the vicinity. Due to cement hydration, there is a~significant temperature increase in the early stages, but it diffuses out quickly. Transport of moisture is much slower, it takes about 30 weeks to reach an~almost steady state. The contours of relative humidity at the final solution step (time one year) are presented in Fig.~\ref{fig:small_dam_paraview}.

Figure~\ref{fig:small_dam}c focuses on the point in the center of the small dam (point $A$ in Fig.~\ref{fig:small_dam}b) and plots time histories of relative humidity, temperature, and cement hydration degree. There is always one graph for the whole simulation time and one for the initial 1/3 of the first week.  

The solutions computed using the \emph{full} and \emph{homogenized} models are almost identical. The \emph{full} model has 25,952 degrees of freedom and runs for 2,873\,s, not considering the preparation time needed to generate the internal structure and the triangulation. The homogenized model has only 360 degrees of freedom and the simulation took around 9\,s. The speed-up factor is larger than 300 (note that different computers and different software were used for each model).

\section{Application \label{sec:application}}

\begin{table}[!b]
	\centering
	\begin{tabular}{rllll}\hline\\[-2.5ex]
	 & $H$ max. & $H$ min.  & $T$ max. & $T$ min.\\\hline\\[-2.5ex]
	water & 1 & 1 & $15^{\circ}$C & $1^{\circ}$C \\
	air & 0.78 & 0.62 & $30^{\circ}$C & $0^{\circ}$C \\
	soil & -- & -- & $19^{\circ}$C & $-1^{\circ}$C \\\hline
	\end{tabular}
	\caption{Extreme values of boundary conditions oscillating during the year according to negative cosine function.\label{tab:BC}}
\end{table}

\begin{figure}[!t]
	\centering
	\includegraphics[width=\textwidth]{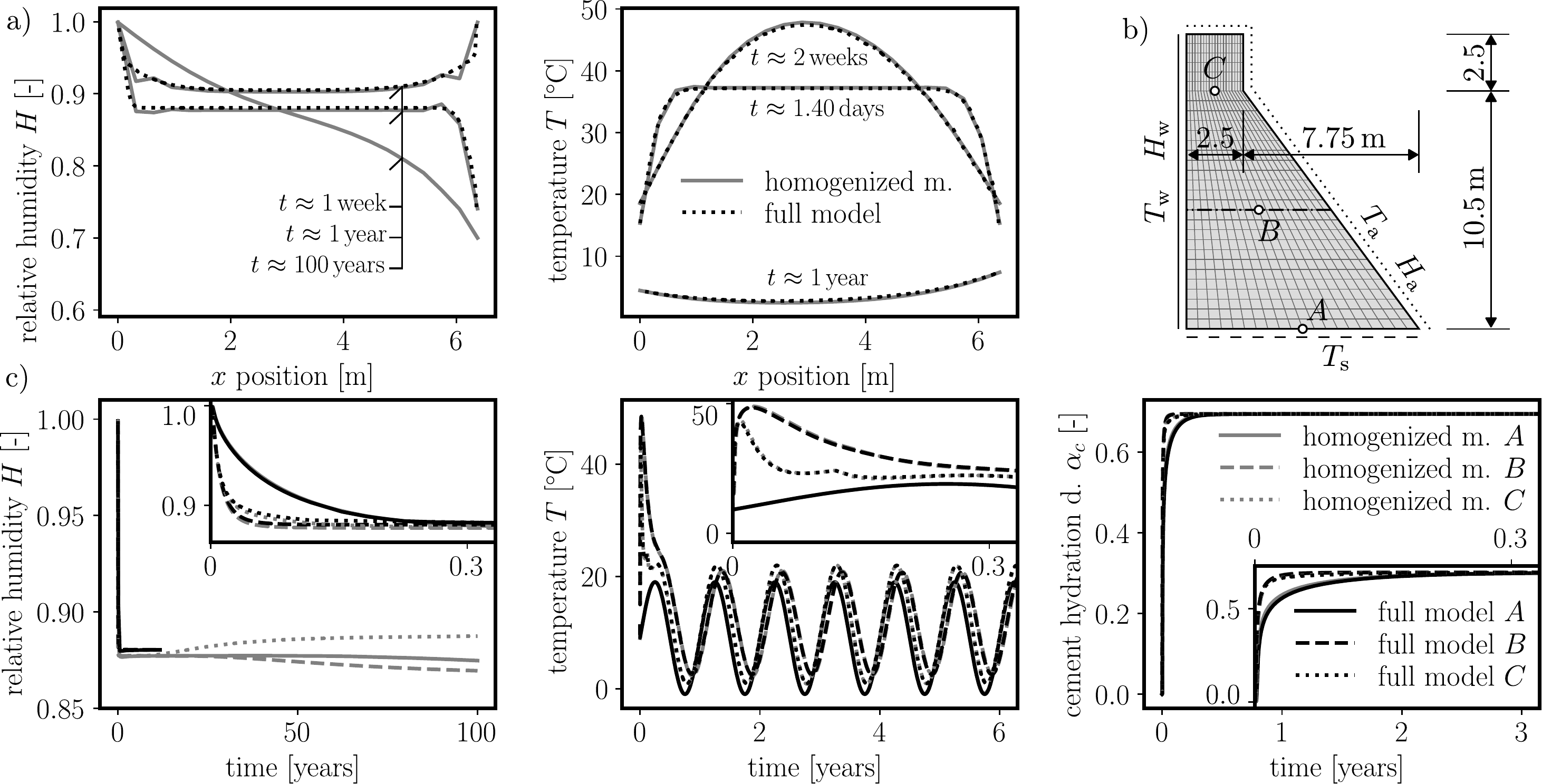}
	\caption{a) Humidity and temperature along the horizontal line through point $B$ at different times; b) dimensions and boundary conditions; c) evolution of humidity, temperature and cement hydration degree at points $A$, $B$ and $C$.  \label{fig:large_dam}}
\end{figure}

The previous examples were meant only as a~verification of the derived homogenization scheme. This section demonstrates that the \emph{homogenized} model is capable to simulate much larger material volumes. The same HTC material model described in~\ref{sec:appendix} is used now to predict the temperature and relative humidity evolution in a~concrete dam with actual dimensions. The geometry and boundary conditions of the dam are adopted from \textcite{gasch2016coupled} and simplified as shown in Fig.~\ref{fig:large_dam}b. The simulated dam has a total width of 10.25\,m, depth of 13\,m and thickness of 1\,m, the maximum aggregate diameter is $d_{\max}=150$\,mm. It is assumed that the dam is kept for 2 weeks under moist-curing conditions. 
After demolding at the end of week 2, it is exposed to ambient relative humidity and temperature for another 4 weeks before filling with the water and starting service at the end of week 6.

The boundary conditions are governed by negative cosine function with a~period of one year and maximum and minimum values reported in Tab.~\ref{tab:BC} for air (a), water (w) and soil (s). These boundary conditions are applied at the right, left and top surfaces according to Fig.~\ref{fig:large_dam}b.  However for the first 2 weeks of curing, $H_{\mathrm{w}}=H_{\mathrm{a}}=1$, then for the next 4 weeks, $H_{\mathrm{w}}$ is considered equal to $H_{\mathrm{a}}$ and finally, after filling the dam, $H_{\mathrm{w}}=1$. Also the water temperature $T_{\mathrm{w}}$ is assumed to be equal to $T_{\mathrm{a}}$ in the first six weeks. Each described change in the boundary conditions is modeled as a one-day linear transition to the next stage. The front and back surfaces of the three-dimensional model have prescribed zero out-of-plane flux for both relative humidity and temperature. Also the out-of-plane flux of humidity at the bottom surface is zero. The initial conditions are $T=15^{\circ}$C and $H=1$. With respect to the negative cosine function, the dam construction begins at 1/4 of the cosine period. The total simulated time is 100 years.

The evolution of relative humidity, temperature and cement hydration degree delivered by the \emph{homogenized} model is shown in Fig.~\ref{fig:large_dam} in gray color. Three different locations $A$, $B$ and $C$ are selected for observation. The plotted time ranges for temperature and cement hydration degree are shortened as the further evolution of these variables follows the pattern from the initial years. The cement hydration degree keeps almost constant for the rest of the simulation time and the temperature keeps oscillating in the same manner as in the first six years. Figure~\ref{fig:large_dam} also shows the relative humidity and temperature profiles along the horizontal line through point $B$ at different time instants. The total computational time of the \emph{homogenized} model was 302\,s using 2,246 degrees of freedom. 

We estimate that the \emph{full} model would need months to compute the same analysis. Instead, we simulated only the first 10 years to verify the \emph{homogenized} model once again. The \emph{full} model had 1,496,652 degrees of freedom, the computational time of the first 10 years was 584,800\,s (about one week). The results are shown in Fig.~\ref{fig:large_dam} in black and they agree well with those obtained by the \emph{homogenized} model. The achieved speed-up factor due to homogenization is roughly four orders of magnitude.

The application demonstrated excellent performance of the proposed homogenization scheme. We admit that similar results might be achieved with tensorial macroscale material description. The effective way of simulating large material volumes with high-fidelity discrete mesoscale models via homogenization becomes very important when the transport or diffusion equation depends on mechanical fields, e.g., crack openings. In such a~case, the mesoscale character of the model is crucial and hardly replaceable by tensorial homogeneous description. Reader is referred to recent paper~\parencite{EliCus22} where homogenization of coupled mechanics and mass transport was developed for the discrete mesoscale model of concrete.

\section{Conclusions}
In this study, a general asymptotic expansion homogenization scheme for mesoscale discrete models of diffusion in heterogeneous materials was derived. The underlying mesoscale model features nonlinearities in all terms of the equation and includes a sink or source term. Based on solid mathematical foundations, the model response is decomposed into continuous macroscopic (slow) and discrete microscopic (fast) components. The slow component is approximated by the finite element method and each integration point contains a~discrete periodic RVE problem that needs to be analyzed in every iteration of the macro model. 

Based on the derivation and the numerical results presented in this paper, the following conclusions can be drawn:
\begin{itemize}
\item The load applied to the RVE problem comes from the pressure gradient at the macroscale as its projection into the normal direction of each RVE discrete conduit element. This projection then serves as an~eigen pressure gradient driving the RVE response.

\item The macroscopic flux vector is computed as the weighted average of oriented fluxes from elements within the RVE. 

\item Both transient and steady-state analyses result in a~steady-state RVE problem. The transient terms appear only at the macroscale. 

\item The RVE problem is linear even when the permeability coefficient depends on the pressure inside the discrete element. The scale separation limits this dependence only to the macroscopic pressure component. The actual degrees of freedom at the RVE level do not affect the permeability coefficients, hence the linearity of the problem. Moreover, for typical nonlinear formulations of the constitutive model where the nonlinearity is multiplicatively added to the basic component, the RVE problem can be pre-computed in advance into an~effective permeability tensor. Subsequently, huge computational time saving is achieved. For the application example considered herein, a speed-up factor of four orders of magnitude was observed.

\item For spatially constant permeability coefficient and a~special type of model geometry (based on Voronoi or power/Laguerre tessellation), the RVE problem disappears, the slow component of pressure becomes the only meaningful part of the solution as expected for homogeneous materials.

\item Verification studies showed an~excellent match between fully resolved and homogenized models for linear and nonlinear as well as steady-state and transient problems. is mainly introduced by poor discretization of the slow solution and can be reduced by adding degrees of freedom to the continuous approximation at the macro scale.

\item The structure of the governing equations is identical to the equations of diffusion, mass transport and conduction of heat or electricity. Therefore, the derived equations can be directly applied to homogenize these problems as well.
\end{itemize}

\section*{Data Availability Statement}
The data that support the findings of this study are available from the corresponding author upon reasonable request.

\section*{CRediT authorship contribution statement}
Jan Eliáš: Methodology, Formal analysis, Software, Validation, Visualization, Writing – original draft, Writing – review \& editing.
Hao Yin: Software, Validation, Writing – review \& editing.
Gianluca Cusatis: Conceptualization, Methodology, Writing – review \& editing.

\section*{Acknowledgements}
Jan Eliáš gratefully acknowledges financial support from the Czech Science Foundation under project no. GA19-12197S.

\appendix
\renewcommand{\thesection}{Appendix \Alph{section}}
\section{Hygro-Thermo-Chemical (HTC) model \label{sec:appendix}}
The HTC model~\parencite{LuzCus09a,LuzCus09b} that simulates moisture and temperature evolution in cementitious materials by considering the evolution of chemical reactions, such as cement hydration and silica fume (pozzolanic) reactions, is briefly reviewed in this section.

Two primary fields, relative humidity $H$ and temperature $T$, and two internal variables, cement hydration degree $\alpha_c$ and silica fume reaction degree $\alpha_s$, are used in the formulation. The mass and heat balance equations read (similarly to Eq.~\ref{eq:balance})
\begin{align}
W\frac{\partial w_{e}}{\partial H} \dot{H}+\sum_{Q\in W} S^{\star} j_H&=Wq_H
&
W\rho c_{t} \dot{T}+\sum_{Q\in W} S^{\star} j_T&=Wq_T
\label{eq:HTCbalance}
\end{align}
where $w_e$ is the evaporable water content, which is a~function of relative humidity $H$, $\alpha_c$ and  $\alpha_s$. $j_H$ is the flux density of water mass per unit time, which is linked to the relative humidity gradient $g_H$ by an~equivalent Darcy's law $j_H=-D_{H}(H_{\lambda}, T_{\lambda}) g_H$. The moisture permeability $D_{H}$ is a~function of  weighted average of relative humidity and temperature $H_{\lambda}$ and $T_{\lambda}$ over the element, respectively. The moisture sink term reads $q_H = \partial w_{e}/\partial \alpha_{c} \dot{\alpha}_{c}+\partial w_{e}/\partial \alpha_{s} \dot{\alpha}_{s}+\partial w_{n}/\partial \alpha_{c} \dot{\alpha}_{c}$, where the non-evaporable water content $w_n$ is a~function of cement hydration degree $\alpha_c$ only, $w_{n}\left(\alpha_{c}\right)=\kappa_{c} \alpha_{c} c$ with $c$ being the cement mass content and $\kappa_{c}=0.253$.
$\rho$ is the density of the material and $c_t$ is the specific heat. $j_T$ is the heat flux density per unit time defined by Fourier's law $j_T=-\kappa g_T $, $\kappa$ is the heat conductivity and $g_T$ is the temperature gradient. The heat source is given by $q_T = \dot{\alpha}_{c} c \tilde{Q}_{c}^{\infty}+\dot{\alpha}_{s} s \tilde{Q}_{s}^{\infty}$, where  $\tilde{Q}_{c}^{\infty}$ is the latent heat of hydration per unit of hydrated mass, $\tilde{Q}_{s}^{\infty}$ is the latent heat of pozzolanic reaction, and $s$ is the silica fume mass content.
The discrete estimation of the temperature or relative humidity gradient is given by Eq.~\eqref{eq:presgrad} using the appropriate primary field. 

The rate of cement hydration degree $\dot{\alpha}_c$ and silica fume reaction degree $\dot{\alpha}_s$ are nonlinear functions of $\alpha_{c}$, $\alpha_s$, $H$ and $T$
\begin{align}
\dot{\alpha}_{c}&=A_{c}\left(\alpha_{c}\right)\beta_{H}\left(H\right) e^{-E_{a c} / (R T)}
&
\dot{\alpha}_{s}&=A_{s}\left(\alpha_{s}\right) e^{-E_{a s} / (R T)}
\end{align}
where $A_{c}\left(\alpha_{c}\right)=A_{c 1}\left(A_{c 2}/\alpha_{c}^{\infty}+\alpha_{c}\right)\left(\alpha_{c}^{\infty}-\alpha_{c}\right) e^{-\eta_{c} \alpha_{c} / \alpha_{c}^{\infty}}$, $\beta_{H}(H)=\left[1+(a-a H)^{b}\right]^{-1}$ and 
$A_{s}\left(\alpha_{s}\right)=A_{s 1}\left(A_{s 2}/\alpha_{s}^{\infty}+\alpha_{s}\right)\left(\alpha_{s}^{\infty}-\alpha_{s}\right) e^{-\eta_{s} \alpha_{s} / \alpha_{s}^{\infty}}$. $E_{ac}/R$, $A_{c 1}$, $A_{c 2}$, $\alpha_{c}^{\infty}$, $\eta_c$, $a$ and $b$ as well as $E_{as}/R$, $A_{s1}$, $A_{s2}$, $\alpha_{s}^{\infty}$ and $\eta_s$ 
are constant material parameters for cement hydration and silica fume reaction, respectively.

The evaporable water content $w_e$ is expressed via a~semi-empirical formula proposed by \textcite{mjornell1997model} to account for the dependencies of evaporable water content to the hydration degree and silica fume reaction degree 
\begin{align}
w_{e}\left(H, \alpha_{c}, \alpha_{s}\right)=& G_{1}\left(\alpha_{c}, \alpha_{s}\right)\left[1-\frac{1}{e^{10\left(g_{1} \alpha_{c}^{\infty}-\alpha_{c}\right) H}}\right] +K_{1}\left(\alpha_{c}, \alpha_{s}\right)\left[e^{10\left(g_{1} \alpha_{c}^{\infty}-\alpha_{c}\right) H}-1\right]
\label{eq:SorptionIsotherm}
\end{align}
where $K_{1}=\left(w_{0}-0.188 \alpha_{c} c+0.22 \alpha_{s} s-G_{1}\left[1-e^{-10\left(g_{1} \alpha_{c}^{\infty}-\alpha_{c}\right)}\right]\right)/\left(e^{10\left(g_{1} \alpha_{c}^{\infty}-\alpha_{c}\right)}-1\right)$, $G_{1}=k_{v g}^{c} \alpha_{c} c+k_{v g}^{s} \alpha_{s} s$. $w_0$, $g_1$, $k_{v g}^{c}$, $k_{v g}^{s}$ are constant material parameters.

The moisture permeability $D_H$ is given by
\begin{align}
    D_{H}(H, T)=\psi(T) D_{1}\left[1+\left(\frac{D_{1}}{D_{0}}-1\right)(1-H)^{n}\right]^{-1}
\end{align}
where $\psi(T)=\exp \left[E_{a d}/(R T_{0})-E_{a d}/(R T)\right]$. $D_{0}$, $D_{1}$, $n$, $E_{a d}/R$ and $T_{0}$ are constant material parameters.

The material parameters used for all calculations with the HTC model in the present paper, taken from~\parencite{LuzCus09b},  are the following: $\rho=\SI{2400}{\kilogram\per\metre\cubed}$, $c_t=\SI{1100}{\joule\per\kelvin\per\kilogram}$, $\kappa=\SI{2.5}{\watt\per\metre\per\kelvin}$, $c=\SI{260.0}{\kilogram\per\metre\cubed}$,  $\tilde{Q}_{c}^{\infty}=\SI{5.2e5}{\joule\per\kilogram}$, $\tilde{Q}_{s}^{\infty}=\SI{7.8e5}{\kilogram\per\metre\cubed}$, $s =\SI{0.0}{\kilogram\per\metre\cubed}$,  $E_{ac}/R=\SI{5490}{\kelvin}$, $A_{c1} =\SI{5.56e4}{\per\second}$, $A_{c2} = \SI{1e-6}{}$, $\alpha_{c}^{\infty} = 0.695$, $\eta_c=6.5$, $a=5.5$, $b=4$, $E_{as}/R=\SI{9620}{\kelvin}$, $A_{s1} = \SI{1.39e10}{\per\second}$, $A_{s2} = 1\times10^{-6}$, $\alpha_{s}^{\infty} = 0.0$, $\eta_s=9.5$, $k_{vg}^c=0.2$, $k_{vg}^s=0.36$, $w_0 = \SI{104}{\kilogram\per\metre\cubed}$, $g_1=1.5$, $\kappa_c=0.253$, $D_0=\SI{6.0e-10}{\kilogram\per\metre\per\second}$, $D_1=\SI{7.0e-8}{\kilogram\per\metre\per\second}$, $n=3$, $E_{ad}/R=\SI{2700}{\kelvin}$, $T_0=\SI{293.15}{\kelvin}$.
They correspond approximately to a~dam concrete mix with cement content $c=\SI{260}{\kilogram\per\metre\cubed}$ and water-to-cement ratio $w/c=0.4$, no silica fume content is involved. A~similar concrete used in drying tests can be found in \textcite{kim1999moisture}.

\printbibliography[heading=bibintoc]

\end{document}